\documentclass[reqno]{amsart}
\usepackage{amssymb,amsmath,amsfonts,amsthm,comment,mathrsfs,times,graphicx}
\usepackage[english]{babel}
\usepackage[all,cmtip]{xy}
\usepackage{lmodern}
\usepackage{enumitem}
\usepackage{geometry}
\newcounter{dummy} \numberwithin{dummy}{section}
\newtheorem{theo}[dummy]{Theorem }
\newtheorem{coro}[dummy]{Corollary}
\newtheorem{lem}[dummy]{Lemma}
\newtheorem{pro}[dummy]{Proposition}
\newtheorem{deft}[dummy]{Definition}
\newtheorem{exe}[dummy]{Example}
\newtheorem{rem}[dummy]{Remark}

\newtheorem{conj}{Conjecture}
\title[ Iwasawa theory of Rubin-Stark units and class group]{Iwasawa theory of Rubin-Stark units and class group}
\author[ Youness Mazigh]{Youness Mazigh}
 \keywords{Iwasawa theory, Euler systems, local conditions.}
\subjclass[2000]{11R23, 11R27, 11R29, 11R42}
\begin{document}
\begin{abstract}
Let $K$ be a totally real number field of degree $r=[K:\mathbb{Q}]$
and let $p$ be an odd rational prime. Let $K_{\infty}$ denote the
cyclotomic $\mathbb{Z}_{p}$-extension of $K$ and let $L_{\infty}$ be
a finite extension of $K_{\infty}$, abelian over $K$. In this
article, we extend results of \cite{Kazim109} relating
characteristic ideal of the $\chi$-quotient of the projective limit
of the ideal class groups to the $\chi$-quotient of the projective
limit of the $r$-th exterior power of units modulo Rubin-Stark
units, in the non semi-simple case, for some
$\overline{\mathbb{Q}_{p}}$-irreductible characters $\chi$ of
$\mathrm{Gal}(L_{\infty}/K_{\infty})$.
\end{abstract}
\maketitle
\section{Introduction}
 Let $K$ be a totally real number field of degree $r=[K:\mathbb{Q}]$. Fix a rational odd
prime $p$ and let $K_{\infty}$ denote the cyclotomic
$\mathbb{Z}_{p}$-extension of $K$; fix  a finite extension
$L_{\infty}$ of $K_{\infty}$, abelian over $K$. Fix also a
decomposition of
$$
\mathrm{Gal}(L_{\infty}/K)=\mathrm{Gal}(L_{\infty}/K_{\infty})\times
\Gamma,\; \Gamma\simeq \mathbb{Z}_{p}.
$$
Then the fields $L:=L_{\infty}^{\Gamma}$ and $K_{\infty}$ are
linearly disjoint over $K$. \vskip 7pt If $F/K$ is a finite abelian
extension of $K$, we write $A(F)$ for the $p$-part of the class
group of $F$, and $\mathcal{E}(F)$ for the group of global units of
$F$. For a $\mathbb{Z}$-module $M$, let $\widehat{M}=\varprojlim
M/p^{n}M$ denote the $p$-adic completion of $M$. Let
$$
A_{\infty}=\varprojlim A(F),\quad\quad
\widehat{\mathcal{E}_{\infty}}=\varprojlim \widehat{\mathcal{E}(F)},
$$
where the projective limit is taken over all finite sub-extensions
of $L_{\infty}$, with respect to  the norm maps. Let
$\Delta=\mathrm{Gal}(L_{\infty}/K_{\infty})$ and let
$$
\xymatrix@=2pc{ \chi : \Delta\ar[r]&
\overline{\mathbb{Q}}_{p}^{\times}}
$$
be a non-trivial $\overline{\mathbb{Q}}_{p}$-irreducible character
of $\Delta$. Let $\mathcal{O}:=\mathbb{Z}_{p}[\chi]$ be the ring
generated by the values of $\chi$ over $\mathbb{Z}_{p}$ and let
$\mathcal{O}(\chi)$ denote the ring $\mathcal{O}$ on which $\Delta$
acts via $\chi$. For any $\mathbb{Z}_{p}[\Delta]$-module $M$, we
define the $\chi$-quotient $M_{\chi}$ by
$$
M_{\chi}:=M\otimes_{\mathbb{Z}_{p}[\Delta]}\mathcal{O}(\chi).
$$
For any profinite group $\mathcal{G}$, we define the Iwasawa algebra
$$
\mathcal{O}[[\mathcal{G}]]:=\varprojlim
\mathcal{O}[\mathcal{G}/\mathcal{H}]
$$
where the projective limit is over all finite quotient
$\mathcal{G}/\mathcal{H}$ of $\mathcal{G}$. In case
$\mathcal{G}=\Gamma$, we shall write
$$
\Lambda:=\mathcal{O}[[\Gamma]].
$$
It is well known that $\Lambda$ is a complete noetherian regular
local domain, with finite residue field of characteristic $p$. The
structure theorem of $\Lambda$-modules shows that for a finitely
generated torsion $\Lambda$-module $M$, there exists an injective
pseudo-isomorphism
$$
\xymatrix@=3pc{ \bigoplus_{i}\Lambda/f_{i}\Lambda\ar@{^{(}->}[r]& M}
$$
with $f_{i}\in \Lambda$, and we define the characteristic ideal of
$M$
$$
\mathrm{char}(M)=\prod_{i}f_{i}\Lambda.
$$
 We denote by $\mathrm{St}_{\infty}:=\varprojlim St_{n}$ the $\mathbb{Z}[[\mathrm{Gal}(L_{\infty}/K)]]$-module
constructed by the Rubin-Stark elements (see Definition
\ref{Definition of Rubin Strak module}). Our objective in this paper
is to compare the characteristic ideal of $(A_{\infty})_{\chi}$ with
the characteristic ideal of
$((\bigwedge^{r}_{\mathbb{Z}_{p}[[\mathrm{Gal(L_{\infty}/K)}]]}\widehat{\mathcal{E}_{\infty}})/\widehat{\mathrm{St}_{\infty}})_{\chi}$.
Let $\Sigma_{\infty}$ be the set of infinite places of $K$ and let
$L_{\chi}$ be the fixed field of $\ker(\chi)$. Let $K(1)$ be the
maximal $p$-extension inside the Hilbert class field of $K$. In the
sequel we will assume  (for simplicity) that
\begin{equation*}\label{simplicity}
L=L_{\chi} \quad\mbox{and}\quad K=L\cap K(1).
\end{equation*}
For a $p$-adic prime $\mathfrak{p}$ of $K$, let
 $\mathrm{Frob}_{\mathfrak{p}}$ denote a Frobenius element at
 $\mathfrak{p}$ inside the absolute Galois group of $K$.
 Assume that
 \begin{center}
  \begin{tabbing}
  \hspace{1cm}\=\hspace{0.4cm} \= \hspace{0.4cm} \= \kill
  \>$(\mathcal{H}_{0})$ \> \> \textit{The number field $L$ is totally real.}
\end{tabbing}
\end{center}
\begin{center}
  \begin{tabbing}
  \hspace{1cm}\=\hspace{0.4cm} \= \hspace{0.4cm} \= \kill
  \>$(\mathcal{H}_{1})$ \> \> \textit{ The extension $L/\mathbb{Q}$ is unramified at $p$.}
\end{tabbing}
\end{center}
\begin{center}
  \begin{tabbing}
  \hspace{1cm}\=\hspace{0.4cm} \= \hspace{0.4cm} \= \kill
  \>$(\mathcal{H}_{2})$ \> \> \textit{ $\chi(\mathrm{Frob}_{\mathfrak{p}})\neq
    1$ for any $p$-adic prime $\mathfrak{p}$ of $K$.}
\end{tabbing}
\end{center}
\begin{center}
  \begin{tabbing}
  \hspace{1cm}\=\hspace{0.4cm} \= \hspace{0.4cm} \= \kill
  \>$(\mathcal{H}_{3})$ \> \> \textit{
   The Leopoldt  conjecture holds for every finite extension $F$ of $L$ in $L_{\infty}$.}
\end{tabbing}
\end{center}
In the semi-simple case B\"{u}y\"{u}kboduk  proved
\begin{theo}(B\"{u}y\"{u}kboduk \cite[Theorem A]{Kazim109})\label{theorem of Kazim}
 Assume that the hypotheses $(\mathcal{H}_{0})-(\mathcal{H}_{3})$
 hold. If $p\nmid [L:K]$, then
$$
\mathrm{char}((A_{\infty})_{\chi})= \mathrm{char}
\bigg(\big((\bigwedge^{r}\widehat{\mathcal{E}_{\infty}})/\widehat{\mathrm{St}_{\infty}}\big)_{\chi}\bigg)
$$
\end{theo}
There is at least  two  ideas behind the proof of such a theorem. On
the one hand, the result may be stated in terms of Selmer groups. On
the other hand, Rubin-Stark elements give rise to Euler systems for
the $p$-adic representation $T=\mathbb{Z}_{p}\otimes
\mathcal{O}(\chi^{-1})$.
 Mazur and Rubin  developed in \cite{MR04} an Euler system and Kolyvagin system
machinery so as to determine the structure of the associated Selmer
groups, in the case where a certain cohomological invariant, called
the Selmer core rank, is one. As an application to Iwasawa theory
B\"{u}y\"{u}kboduk obtains a divisibility relation between the
characteristic ideals of the projective limit of these Selmer
groups, which transforms into equality if the corresponding
Kolyvagin system is primitive. For more detail, see \cite{Kazim107}.
He then applied this theory to the proof of Theorem \ref{theorem of
Kazim} by  constructing a primitive Kolyvagin system from
Rubin-Stark elements and  Selmer groups of core rank one.
\begin{rem}
 Mazur and Rubin introduced in \cite{MR 16} the notion of Stark
system/Kolyvagin system of rank $r$. They used this notion to
determine the structure of Selmer groups, when the core rank is
greater than one.
\end{rem}
  \vskip 6pt
 In this paper we prove
\begin{theo}\label{AMO2}
 Assume that the hypotheses $(\mathcal{H}_{0})-(\mathcal{H}_{3})$
 hold. Then
$$
\mathrm{char}((A_{\infty})_{\chi})\quad \mbox{divides}\quad p^{d}
\mathrm{char}
\big(((\bigwedge^{r}\widehat{\mathcal{E}_{\infty}})/\widehat{\mathrm{St}_{\infty}})_{\chi}\big)
$$
where $d=\displaystyle{\max_{\mathfrak{p}\mid
p}}\{v_{p}(1-\chi(\mathrm{Frob}_{\mathfrak{p}}))\}$.
\end{theo}
We take our inspiration from \cite{Kazim109}. But if $p\mid [L:K]$
the results of \cite{MR04}, \cite{MR 16} and \cite{Kazim107} do not
apply, since the notion of core rank is not defined. Therefore, we
are led to use the theory of Euler systems  exposed in
\cite{Rubin00}. In particular, we construct an ad-hoc Selmer
structure and an associated  Kolyvagin system. We have to use the
structure of the semi-local units, \;cf.\;Theorem \ref{Pro H1IW rank
}. But this is already known, thanks to Greither who applied
Coleman's theory in \cite{Greither96} to determine this structure.
In the semi-simple case B\"{u}y\"{u}kboduk
 used a weak structure theorem of  Colmez-Cherbonnier, obtained by
 using the theory of $(\varphi,\Gamma)$-modules.
\section{\bf Selmer structures}
In this section we recall  some definitions concerning the notion of
 Selmer structure introduced by Mazur and Rubin in \cite{MR04} and
 \cite{MR 16}. For any field $k$  and a fixed separable algebraic closure
$\overline{k}$ of $k$, we write
$G_{k}:=\mathrm{Gal}(\overline{k}/k)$ for the Galois group of
$\overline{k}/k$. Let $\mathcal{O}$ be the ring of integers of a
finite extension $\Phi$ of $\mathbb{Q}_{p}$ and  let $D$ denote the
divisible module $\Phi/\mathcal{O}$. For a $p$-adic representation
$T$ with coefficients in $\mathcal{O}$, we  define
$$ D(1)=D\otimes \mathbb{Z}_{p}(1),\quad\quad
T^{\ast}=\mathrm{Hom}_{\mathcal{O}}(T,D(1)),
$$
where $\mathbb{Z}_{p}(1):=\varprojlim \mu_{p^{n}}$ is the Tate
module.\vskip 6pt Let $F$ be a number field, and for a place $w$ of
$F$, let $F_{w}$ denote the completion of $F$ at the  place $w$. Let
us recall the local duality theorem c.f.\cite[Corollary I.2.3
]{Milne}: For $i=0,1,2$, there is a perfect  pairing
\begin{equation}\label{local dualty theorem}
\begin{array}{cccccccc}
  H^{2-i}(F_{w},T) & \times & H^{i}(F_{w},T^{\ast})& \xymatrix@=3pc{\ar[r]^-{\langle\;,\;\rangle_{w}}&} &H^{2}(F_{w},D(1))\cong D,\mbox{if $v$ is finite},\\
 &  &   &  &\\
  \widehat{H}^{2-i} (F_{w},T)& \times & \widehat{H}^{i}(F_{w},T^{\ast}) & \xymatrix@=3pc{\ar[r]^-{\langle\;,\;\rangle_{w}}&} &
\widehat{H}^{2}(F_{w},D(1)),\mbox{if $v$ is infinite}
\end{array}
\end{equation}
where $\widehat{H}^{\ast}(F_{w},.)$ denotes the Tate cohomology
group.\vskip 6pt
\begin{deft} Let $T$ be a $p$-adic representation of $G_{F}$ with
coefficients in $\mathcal{O}$ and let $w$ be a non $p$-adic prime of
$F$. A local condition $\mathcal{F}$ at the prime $w$  on $T$ is a
choice of an $\mathcal{O}$-module $H^{1}_{\mathcal{F}}(F_{w},T)$ of
$H^{1}(F_{w},T)$. For the $p$-adic primes, a local condition at $p$
will be a choice of an $\mathcal{O}$-submodule
$H^{1}_{\mathcal{F}}(F_{p},T)$ of the semi-local cohomology group
$$
H^{1}(F_{p},T):=\oplus_{w\mid p}H^{1}(F_{w},T).
$$
\end{deft}
Let $I_{w}$ denote the inertia subgroup of $G_{F_{w}}$. We say that
$T$ is unramified at $w$ if the inertia subgroup $I_{w}$ of $w$ acts
trivially on $T$. We assume in the sequel that $T$ is unramified
outside a finite set of places of $F$.
\begin{deft}
A Selmer structure $\mathcal{F}$ on $T$ is a collection of the
following data:
\begin{itemize}
    \item a finite set $\Sigma(\mathcal{F})$ of $F,$ including all infinite places and primes above $p$, and all primes
where $T$ is ramified;
    \item for every $w\in \Sigma(\mathcal{F})$, a local condition on
    $T$.
\end{itemize}
\end{deft}
 If $\mathcal{F}$ is a Selmer structure on $T$, we define the Selmer group
$H^{1}_{\mathcal{F}}(F,T)\subset H^{1}(F,T)$ to be the kernel of the
sum of the restriction maps
\begin{equation}\label{Selmer group definition}
\xymatrix@=2pc{H^{1}(G_{\Sigma(\mathcal{F})}(F),T)\ar[r]&
\displaystyle{{\bigoplus}_{w\in\Sigma(\mathcal{F})}}(H^{1}(F_{w},T)/H^{1}_{\mathcal{F}}(F_{w},T))},
\end{equation}
where
$G_{\Sigma(\mathcal{F})}(F):=\mathrm{Gal}(F_{\Sigma(\mathcal{F})}/F)$
is the Galois group of the maximal algebraic extension of $F$ which
is unramified outside $\Sigma(\mathcal{F})$.\vskip 6pt A  Selmer
structure $\mathcal{F}$ on $T$ determines a Selmer structure
$\mathcal{F}^{\ast}$ on $T^{\ast}$. Namely,

$$ \Sigma(\mathcal{F})=\Sigma(\mathcal{F}^{\ast}),\quad
H^{1}_{\mathcal{F}^{\ast}}(F_{w},T^{\ast}):=H^{1}_{\mathcal{F}}(F_{w},T)^{\perp},\;\;
\mbox{if $w\in\Sigma(\mathcal{F}^{\ast})-\Sigma_{p}$}
$$
 under the local Tate pairing $\langle\;,\;\rangle_{w}$
and
$$
H^{1}_{\mathcal{F}^{\ast}}(F_{p},T^{\ast}):=H^{1}_{\mathcal{F}}(F_{p},T)^{\perp},
$$
 under the pairing $\sum_{w\mid p}\langle\;,\;\rangle_{w}$.\vskip 6pt
\begin{exe}\label{Example unramified local condition}
 Let $w$ be a place of $F$ and  let $F_{w}^{ur}$ denote the maximal
 unramified extension of $F_{w}$. Define the subgroup of universal norm
$$
H^{1}(F_{w},T)^{u}=\bigcap_{F_{w}\subset k\subset
F_{w}^{ur}}\mathrm{cor}_{k,F_{w}}H^{1}(k,T),
$$
 where the intersection is over all finite unramified extensions $k$ of $F_{w}$.
Let $H^{1}(F_{w},T)^{u,sat}$ denote the $\mathcal{O}$-saturation of
$H^{1}(F_{w},T)^{u}$ in $H^{1}(F_{w},T)$, $\mathrm{i.e}$,
$H^{1}(F_{w},T)/H^{1}_{\mathcal{F}_{ur}}(F_{w},T)$ is a free
$\mathcal{O}$-module and
$H^{1}_{\mathcal{F}_{ur}}(F_{w},T)/H^{1}(F_{w},T)^{u}$ has finite
length. Following \cite[Defintition 5.1]{MR 16}, we define the
unramified local condition $\mathcal{F}_{ur}$
 by
\begin{equation}\label{Defintion of unramified condition}
H^{1}_{\mathcal{F}_{ur}}(F_{w},T)=H^{1}(F_{w},T)^{u,sat},\;\mbox{if
$w\nmid p$},\quad\mbox{and}\quad
H^{1}_{\mathcal{F}_{ur}}(F_{p},T)={\bigoplus}_{\mathfrak{p}\mid
p}H^{1}(F_{\mathfrak{p}},T)^{u,sat}
\end{equation}
For any $G_{F_{w}}$-module $M$, we define the subgroup of unramified
cohomology classes $H^{1}_{ur}(F_{w},M)\subset H^{1}(F_{w},M)$ by

$$
H^{1}_{ur}(F_{w},M)=\ker(\xymatrix@=2pc{H^{1}(F_{w},M)\ar[r]&
H^{1}(I_{w},M)}).
$$
For future use, we record here the following well-known properties
of unramified cohomology
\begin{enumerate}[label=(\roman*)]
    \item  $H^{1}_{\mathcal{F}_{ur}^{\ast}}(F_{w},T^{\ast})=H^{1}_{ur}(F_{w},T^{\ast})_{div},\quad H^{1}_{\mathcal{F}_{ur}^{\ast}}(F_{p},T^{\ast})={\bigoplus}_{\mathfrak{p}\mid
p}H^{1}_{ur}(F_{\mathfrak{p}},T^{\ast})_{div}.$
    \item If $w\nmid p$ and $T$ is
unramified at $w$, then
$$
H^{1}_{\mathcal{F}_{ur}}(F_{w},T)=H^{1}_{ur}(F_{w},T)\quad\mbox{and}\quad
H^{1}_{\mathcal{F}_{ur}^{\ast}}(F_{w},T^{\ast})=H^{1}_{ur}(F_{w},T^{\ast}).
$$
\end{enumerate}
 where for an
abelian group $A$, $A_{div}$ denotes the maximal divisible subgroup
of $A$.\\
The assertion $(\textrm{i})$ follows from \cite[\S 2.1.1,
Lemme]{PR92} and the assertion $(\textrm{ii})$ can be deduced from
\cite[Lemma 1.3.5]{Rubin00}.
\end{exe}
\section{\bf Iwasawa Theory}
 Fix a totally real number field $K$. Let $r=[K:\mathbb{Q}]$ and
$K_{\infty}=\bigcup_{n\geq 0}K_{n}$ denote the cyclotomic
$\mathbb{Z}_{p}$-extension of $K$. Assume that all algebraic
extensions of $K$ are contained in a fixed  algebraic closure
$\overline{\mathbb{Q}}$ of $\mathbb{Q}$. If $F$ is a finite
extension of $K$ and $w$ is a place of $F$, fix a place
$\overline{w}$ of $\overline{\mathbb{Q}}$ lying above $w$. The
decomposition (resp. inertia) group of $\overline{w}$ in
$\overline{\mathbb{Q}}/F$ is denoted  by $D_{w}$ (resp. $I_{w}$). If
$v$ is a prime of $K$ and $F$ is a Galois extension of $K$, we
denote the decomposition group of $v$ in $F/K$ by $D_{v}(F/K)$.
Recall that
$$
\xymatrix@=2pc{ \chi : G_{K}\ar[r]& \mathcal{O}^{\times}}
$$
is a non-trivial $\overline{\mathbb{Q}}_{p}$-irreducible character,
factoring through a finite abelian  extension $L$ of $K$. Assume
that $L$ and $K_{\infty}$ are linearly disjoint over $K$. Let
$L_{n}=LK_{n}$ and let $L_{\infty}=LK_{\infty}$ be the cyclotomic
$\mathbb{Z}_{p}$-extension of $L$. In the sequel, we will denote by
$T$ the $p$-adic representation
$$
T=\mathbb{Z}_{p}(1)\otimes \mathcal{O}(\chi^{-1}).
$$
Let $\Sigma$ be a finite set of places of $K$ containing all
infinite places, all $p$-adic places and all places where $T$ is
ramified. If $F$ is an extension of $K$, we denote also by $\Sigma$
the set of places of $F$ lying above places in $\Sigma$.
\begin{deft}\label{definition of canonical and strict conditions}
Let $F$ be a finite extension of $K$ and let $\Sigma_{p}$ denote the
set of $p$-adic primes of $K$. Following \cite[Example 5.1]{MR 16},
we define the canonical
 and the strict Selmer structures $\mathcal{F}_{can}$ and
$\mathcal{F}_{str}$ on $T$, by
\begin{itemize}
    \item
    $\Sigma(\mathcal{F}_{can})=\Sigma(\mathcal{F}_{str})=\Sigma$.
    \item if $w\nmid p$,
    $H^{1}_{\mathcal{F}_{can}}(F_{w},T)=H^{1}_{\mathcal{F}_{str}}(F_{w},T):=H^{1}_{\mathcal{F}_{ur}}(F_{w},T)$.
    \item $H^{1}_{\mathcal{F}_{can}}(F_{p},T):=H^{1}(F_{p},T)$, and
    $H^{1}_{\mathcal{F}_{str}}(F_{p},T)=0$.
\end{itemize}
where $\mathcal{F}_{ur}$ is the unramified local condition defined
in $(\ref{Defintion of unramified condition})$.
\end{deft}
Let $F^{\prime}/F$ be a finite extension of $F$. Remark that for
$\mathcal{F}=\mathcal{F}_{can}$, $\mathcal{F}=\mathcal{F}_{ur}$ or
$\mathcal{F}=\mathcal{F}_{str}$, we have
$$
\mathrm{cor}_{F^{\prime}_{w^{\prime}},F_{w}}(H^{1}_{\mathcal{F}}(F^{\prime}_{w^{\prime}},T))\subset
H^{1}_{\mathcal{F}}(F_{w},T)\quad \mbox{and}\quad
\mathrm{res}_{F_{w},F^{\prime}_{w^{\prime}}}(H^{1}_{\mathcal{F}}(F_{w},T^{\ast}))\subset
H^{1}_{\mathcal{F}}(F^{\prime}_{w^{\prime}},T^{\ast})
$$
where $w^{\prime}\mid w$ and
$\mathrm{cor}_{F^{\prime}_{w^{\prime}},F_{w}}$ (resp.
$\mathrm{res}_{F_{w},F^{\prime}_{w^{\prime}}}$) denote the
corestriction (resp.  restriction) map. For these  local conditions,
we write
$$
H^{1}_{\mathcal{F}}(FK_{\infty},T):=\varprojlim_{n}H^{1}_{\mathcal{F}}(FK_{n},T),\quad
H^{1}_{\mathcal{F}^{\ast}}(FK_{\infty},T^{\ast}):=\varinjlim_{n}H^{1}_{\mathcal{F}^{\ast}}(FK_{n},T^{\ast})
$$
 where the projective (resp. injective) limit is taken with respect to
 the corestriction (resp. restriction) maps.\vskip 6pt
 Recall that for any Galois extension $F/F^{\prime}$ of number
 fields, $K\subset F^{\prime}\subset F$, the restriction map induces an isomorphism
 \begin{equation}\label{lemme non semi simple}
\xymatrix@=2pc{ \mathrm{res}: H^{1}(F^{\prime},T)\ar[r]^-{\sim}&
H^{1}(F,T)^{\mathrm{Gal}(F/F^{\prime})}},
 \end{equation}
 c.f.\,\cite[Lemme 4.3]{AMO1}.
\begin{lem}\label{Leopoldt}
 We have
$$
H^{1}_{\mathcal{F}_{str}}(K_{\infty},T)=0.
$$
\end{lem}
\noindent \textbf{Proof.} This equality is implicit in
\cite[Proposition 2.12]{Kazim109}. It is a consequence of the weak
Leopoldt conjecture. Indeed, using $(\ref{lemme non semi simple})$
and passing to the inverse limit, we obtain
$$
\varprojlim_{n} H^{1}(G_{\Sigma}(K_{n}),T)\cong(\varprojlim_{n}
H^{1}(G_{\Sigma}(L_{n}),T))^{\mathrm{Gal}(L_{\infty}/K_{\infty})},
$$
then $\xymatrix@=2pc{
H^{1}_{\mathcal{F}_{str}}(K_{\infty},T)\ar@{^{(}->}[r]&
H^{1}_{\mathcal{F}_{str}}(L_{\infty},T)^{\mathrm{Gal}(L_{\infty}/K_{\infty})}}$.
This shows that it suffices to prove
$H^{1}_{\mathcal{F}_{str}}(L_{\infty},T)=0$. Indeed, for
$w\in\Sigma-\Sigma_{p}$, Proposition $\mathrm{B}.3.2$ of
\cite{Rubin00} and $(i)$ of  example \ref{Example unramified local
condition} show that
$$\displaystyle{\varprojlim_{n}}H^{1}(L_{n,w},T)\cong
\displaystyle{\varprojlim_{n}}H^{1}_{ur}(L_{n,w},T)
\quad\mbox{and}\quad
H^{1}_{\mathcal{F}_{ur}}(L_{n,w},T)=H^{1}_{ur}(L_{n,w},T).
$$
Then by definition of $H^{1}_{\mathcal{F}_{str}}(L_{\infty},T)$, we
have an exact sequence
$$
\xymatrix@=2pc{0\ar[r]&H^{1}_{\mathcal{F}_{str}}(L_{\infty}),T)\ar[r]&\varprojlim_{n}H^{1}(G_{\Sigma}(L_{n},T)\ar[r]&
\varprojlim_{n}H^{1}(L_{n,p},T)}.
$$
Since $L_{\infty}$ is the cyclotomic $\mathbb{Z}_{p}$-extension of
$L$, then the weak Leopold conjecture is true for $L_{\infty}/L$,
e.g.\,\cite[Theorem 10.3.25]{NSW91}. Therefore, using
$\chi(G_{L})=1$, we get
$$
\varprojlim_{n}H^{1}_{\mathcal{F}_{str}}(L_{n},T)=0.
$$
 \hfill
$\square$\vskip 6pt Let $n$ be a nonnegative integer, we write
$A_{n}$ for the $p$-part of the class group of $L_{n}$ and
$\mathcal{E}^{\prime}_{n}$ for the group of $p$-units of $L_{n}$.
Let
$$
A_{\infty}:=\varprojlim_{n}A_{n}\quad \mbox{and}\quad
\widehat{\mathcal{E}^{\prime}_{\infty}}:=\varprojlim_{n}\widehat{\mathcal{E}^{\prime}}_{n},
$$
 where all inverse limits are taken with respect to norm maps. It is well known that
$$
\varprojlim_{n}H^{1}(G_{\Sigma_{p\infty}}(L_{n}),\mathbb{Z}_{p}(1))\cong
\widehat{\mathcal{E}^{\prime}_{\infty}}
$$
Then, by $(\ref{lemme non semi simple})$, we have
\begin{equation}\label{Fcan and units }
H^{1}_{\mathcal{F}_{can}}(K_{\infty},T)\cong
\varprojlim_{n}H^{1}(K_{n},T)\cong
(\widehat{\mathcal{E}^{\prime}_{\infty}}\otimes_{\mathbb{Z}_{p}}\mathcal{O}(\chi^{-1}))^{\mathrm{Gal}(L_{\infty}/K_{\infty})}.
\end{equation}
\begin{pro}\label{Proposition fcan is free of rank r}
Suppose that every infinite place of $K$ is completely decomposed in
$L/K$. Then the $\Lambda$-module
$H^{1}_{\mathcal{F}_{can}}(K_{\infty},T)$ is free of rank
$r=[K:\mathbb{Q}].$
\end{pro}
\noindent \textbf{Proof.} Using Dirichet's unit theorem and the
assumption $(\mathcal{H}_{0})$, we get
$$
\mathrm{rank}_{\mathcal{O}}((\widehat{\mathcal{E}^{\prime}_{n}}
\otimes_{\mathbb{Z}_{p}}\mathcal{O}(\chi^{-1}))^{\mathrm{Gal}(L_{n}/K_{n})})=r.p^{n}+t
$$
where $t$ is a nonnegative integer independent of $n$. Then by
\cite[Theorem]{Greither94}, we see that
$$
\mathrm{rk}_{\Lambda}((\widehat{\mathcal{E}^{\prime}_{\infty}}
\otimes_{\mathbb{Z}_{p}}\mathcal{O}(\chi^{-1}))^{\mathrm{Gal}(L_{\infty}/K_{\infty})})=r.
$$
This finishes the proof of the proposition.\hfill $\square$\vskip
6pt

\begin{pro}\label{proposition class group and Fur}
The $\mathcal{O}[\mathrm{Gal}(L_{n}/K)]$-modules
$H^{1}_{\mathcal{F}_{ur}^{\ast}}(L_{n},T^{\ast})$ and
$\mathrm{Hom}(A_{n},T^{\ast})$ are isomorphic.
\end{pro}
\noindent \textbf{Proof.} See $§ 6.1$ of \cite{MR04}.\hfill
$\square$\vskip 5pt
%
%
 For an $\mathcal{O}$-module $M$, we denote by $M^{\vee}:=\mathrm{Hom}_{\mathcal{O}}(M,D)\cong
  \mathrm{Hom}_{\mathbb{Z}_{p}}(M,\mathbb{Q}_{p}/\mathbb{Z}_{p})$ its
  Pontryagin dual.
\begin{pro}\label{Proposition Fcan pseudo isomo to co inv}
The $\Lambda$-modules
$H^{1}_{\mathcal{F}^{\ast}_{can}}(K_{\infty},T^{\ast})^{\vee}$ and
$(H^{1}_{\mathcal{F}^{\ast}_{can}}(L_{\infty},T^{\ast})^{\vee})_{\mathrm{Gal}(L_{\infty}/K_{\infty})}$
are pseudo-isomorphic.
\end{pro}
\noindent \textbf{Proof.} This is \cite[Proposition
3.8]{AMO1}.\hfill $\square$
\begin{lem}\label{chi quotient of class group and canonique
conditio} If one of the hypotheses $(\mathcal{H}_{2})$ or
$(\mathcal{H}_{3})$ holds then
$$
\mathrm{char}((A_{\infty})_{\chi})\quad\mbox{divides}\quad
\mathrm{char}(H^{1}_{\mathcal{F}^{\ast}_{can}}(K_{\infty},T^{\ast})^{\vee}).
$$
\end{lem}
\noindent \textbf{Proof.} Consider the exact sequence
$$
\xymatrix@=2pc{H^{1}_{\mathcal{F}_{ur}^{\ast}}(L_{n,p},T^{\ast})^{\vee}\ar[r]&
H^{1}_{\mathcal{F}_{ur}^{\ast}}(L_{n},T^{\ast})^{\vee}\ar[r]&
H^{1}_{\mathcal{F}_{can}^{\ast}}(L_{n},T^{\ast})^{\vee}\ar[r]&0}.
$$
Since
$$
H^{1}_{\mathcal{F}_{ur}^{\ast}}(L_{n,p},T^{\ast})\cong
\bigoplus_{w\mid p}\mathrm{Hom}(D_{w}/I_{w},T^{\ast})
$$
 Then
$$
H^{1}_{\mathcal{F}_{ur}^{\ast}}(L_{n,p},T^{\ast})^{\vee}\cong
\bigoplus_{v\mid
p}\mathcal{O}(\chi^{-1})[\mathrm{Gal}(L_{n}/K)/D_{v}(L_{n}/K)].
$$
Passing to the projective limit and taking the
$\Delta$-co-invariants, we get
$$
(\mathcal{O}(\chi^{-1})[\mathcal{G}/D_{v}(L_{\infty}/K)])_{\Delta}\simeq\left\{
                                                                          \begin{array}{ll}
                                                                            \mbox{finite}, & \hbox{if $\chi(D_{v}(L/K))\neq 1$;} \\
                                                                            \mathcal{O}[\mathrm{Gal}(K_{\infty}/K)/D_{v}(K_{\infty}/K)],
                                                                             & \hbox{if  $\chi(D_{v}(L/K))=1$.}
                                                                          \end{array}
                                                                        \right.
$$
where $\Delta=\mathrm{Gal}(L_{\infty}/K_{\infty})$. Using
Proposition \ref{proposition class group and Fur} and Proposition
\ref{Proposition Fcan pseudo isomo to co inv}, we obtain that
$$
\mathrm{char}((A_{\infty})_{\chi})\quad\mbox{divides}\quad
\mathcal{J}^{s}\mathrm{char}(H^{1}_{\mathcal{F}^{\ast}_{can}}(K_{\infty},T^{\ast})^{\vee})
$$
where $\mathcal{J}$ is the augmentation ideal of $\Lambda$ and
$s=\#\{ v\mid p \;;\; \chi(\mathrm{Frob}_{v})=1  \}$. It is well
known that $\mathrm{char}((A_{\infty})_{\chi})$ is prime to
$\mathcal{J}$ in the cyclotomic $\mathbb{Z}_{p}$-extension if
$(\mathcal{H}_{3})$ is satisfied. This concludes the proof of the
lemma. \hfill $\square$\vskip 6pt
\section{\bf Euler systems of Rubin-Stark units}
In this section, we construct an Euler system in the sense of
\cite[Definition 2.1.1]{Rubin00} for the $p$-adic representation
$\mathbb{Z}_{p}(1)\otimes\mathcal{O}(\chi^{-1})$, coming from the
elements predicted by Rubin-Stark conjecture \cite[Conjecture
$B^{\prime}$]{Rubin96}.\vskip 6pt We set some notation. Let $K$ be a
number field and let $F$ be a finite abelian extension of $K$. Fix a
finite set $S$ of places of $K$ containing all infinite places and
all places ramified in $F/K$, and a second finite set $\mathcal{T}$
of places of $K$, disjoint from $S$. Let $G=\mathrm{Gal}(F/K)$ and
$\widehat{G}=\mathrm{Hom}(G,\mathbb{C}^{\times})$. If $\rho\in
\widehat{G}$ we define the modified Artin $L$-function attached to
$\rho$ by
$$
L_{S,\mathcal{T}}(s,\rho)=\prod_{\mathfrak{p}\not\in
S}(1-\rho(\mathrm{Frob}_{\mathfrak{p}})\mathbf{N}\mathfrak{p}^{-s})^{-1}
\prod_{\mathfrak{p}\in
\mathcal{T}}(1-\rho(\mathrm{Frob}_{\mathfrak{p}})\mathbf{N}\mathfrak{p}^{1-s})
$$
where $\mathrm{Frob}_{\mathfrak{p}}\in G$ is the Frobenius of the
(unramified) prime $\mathfrak{p}$.\par For each $\rho\in
\widehat{G}$, there is an idempotent
$$
e_{\rho}=\frac{1}{|G|}\sum_{\sigma\in G}\rho(\sigma)\sigma^{-1}\in
\mathbb{C}[G].
$$
Following \cite{Tate84} we define the Sticklberger element
$$
\Theta_{S,\mathcal{T}}(s)=\Theta_{S,\mathcal{T},F/K}(s)=\sum_{\rho\in\widehat{G}}L_{S,\mathcal{T}}(s,\rho^{-1})e_{\rho}
$$
which we view as a $\mathbb{C}[G]$-valued meromorphic function on
$\mathbb{C}$.
  Let $\rho\in \widehat{G}$ and let
$r_{S}(\rho)$ be  the order of vanishing of
$L_{S,\mathcal{T}}(s,\rho)$ at $s=0$. Recall that
$$
r_{S}(\rho)=\mathrm{ord}_{s=0}L_{S,\mathcal{T}}(s,\rho)=\left\{
                                              \begin{array}{ll}
                                                |\{v\in S\;:\; \rho(D_{v}(F/K))=1\}|, & \hbox{$\rho\neq 1$;} \\
                                                |S|-1, & \hbox{$\rho=1$.}
                                              \end{array}
                                            \right.
$$
(see for example \cite{Tate84} Proposition $\textrm{I}.3.4$), where
$D_{v}(F/K)$ is the decomposition group of $v$ relative to
$F/K$.\vskip 6pt
 Before stating the Rubin-Stark conjecture we record some hypotheses
 $\mathbf{H}(F/K,S,\mathcal{T},r)$:
\begin{enumerate}
    \item $S$ contains all the infinite primes of $K$ and all the
    primes which ramify in $F/K$;
    \item $S$ contains at least $r$ places which split completely in
    $F/K$;
    \item $|S|\geq r+1$;
    \item $\mathcal{T}\neq \emptyset$, $S\cap \mathcal{T}=\emptyset$ and $U_{S,\mathcal{T}}(F)$ is
    torsion-free,
\end{enumerate}
here $U_{S,\mathcal{T}}(F)$ is the group of $S$-units of $F$ which
are congruent to $1$ modulo all the primes in $\mathcal{T}$.
\begin{rem}
Conditions $(2)$ and $(3)$ ensure that
$s^{-r}\Theta_{S,\mathcal{T}}(s)$ is holomorphic at $s=0$. Condition
$(4)$ is easily satisfied. For example, if $\mathcal{T}$ contains
primes of two different residue characteristics.
\end{rem}
We will identify
$\mathrm{Hom}_{\mathbb{Z}[G]}(U_{S,\mathcal{T}}(F),\mathbb{Z}[G])$
with a submodule of $\mathrm{Hom}_{\mathbb{C}[G]}(\mathbb{C}\otimes
U_{S,\mathcal{T}}(F), \mathbb{C}[G])$. For any  $r$-tuple
$(\phi_{1},\cdots,\phi_{r})\in\mathrm{Hom}_{\mathbb{Z}[G]}(U_{S,\mathcal{T}}(F),\mathbb{Z}[G])^{r}$,
we define a $\mathbb{C}[G]$-morphism
$$
\xymatrix@=3.5pc{
\mathbb{C}\otimes\bigwedge^{r}_{\mathbb{Z}[G]}U_{S,\mathcal{T}}(F)
\ar[r]^-{\phi_{1}\wedge\cdots\wedge \phi_{r}}& \mathbb{C}[G]}
$$
by
$$
\phi_{1}\wedge\cdots\wedge \phi_{r}(u_{1}\wedge\cdots\wedge
u_{r}):=\displaystyle{\det_{1\leq i,j\leq r}}(\phi_{i}(u_{j})),
$$
for any $u_{1},\cdots, u_{r}\in U_{S,\mathcal{T}}(F)$.\vskip 6pt
 For any $\mathbb{Z}[G]$-module $M$ with trivial
$\mathbb{Z}$-torsion and any positive integer $r$, we let
$$
M_{r,S}:=\big\{ x\in M\;|\;\mbox{ $e_{\rho}.x=0$ in
$\mathbb{C}\otimes M$ for all $\rho\in \widehat{G}$ such that
$r_{S}(\rho)> r$}\big\}.
$$
Assuming that $(F/K,S,\mathcal{T},r)$ satisfies hypotheses
$\mathbf{H}(F/K,S,\mathcal{T},r)$ and $r\geq 1$, let
$\Lambda_{S,\mathcal{T}}$ the $\mathbb{Z}[G]$-submodule of
$\mathbb{Q}\otimes\bigwedge^{r}_{\mathbb{Z}[G]}U_{S,\mathcal{T}}(F)$
defined by
 $$ \Lambda_{S,\mathcal{T}}:={\bigg\{\substack{ x\in
(\mathbb{Q}\otimes\bigwedge^{r}_{\mathbb{Z}[G]}U_{S,\mathcal{T}}(F))_{r,S}\quad|\quad
(\phi_{1}\wedge\cdots\wedge \phi_{r})(x)\in \mathbb{Z}[G],\\  \\
\forall
\phi_{1},\cdots,\phi_{r}\in\mathrm{Hom}_{\mathbb{Z}[G]}(U_{S,\mathcal{T}}(F),\mathbb{Z}[G])
}\bigg\}} .$$
\begin{rem}\label{Remark Rubin lattice}
It is immediate that for $r=1$, we have
$\Lambda_{S,\mathcal{T}}=(\widetilde{U_{S,\mathcal{T}}(F)})_{1,S}$,
where for a $\mathbb{Z}[G]$-module $M$, $\widetilde{M}$ denotes the
image of $M$ via the canonical morphism $M\longrightarrow
\mathbb{Q}\otimes M$. For a general $r\geq 1$, we have inclusions
$$
\mid G\mid^{n}\Lambda_{S,\mathcal{T}}\subset
(\widetilde{\bigwedge^{r}_{\mathbb{Z}[G]}U_{S,\mathcal{T}}})_{r,S}\subset
\Lambda_{S,\mathcal{T}},
$$
for sufficiently large positive integer $n$. Since
$U_{S,\mathcal{T}}(F)$ has finite index in $U_{S}(F)$,
$$
\mathbb{Q}\otimes \Lambda_{S,\mathcal{T}}= (\mathbb{Q}\otimes
\bigwedge^{r}_{\mathbb{Z}[G]}U_{S,\mathcal{T}}(F))_{r,S}=(\mathbb{Q}\otimes
\bigwedge^{r}_{\mathbb{Z}[G]}U_{S}(F))_{r,S}.
$$
Remark also that the module  $(\mathbb{Q}U_{S}(F))_{r,S}$ is
isomorphic to $(\mathbb{Q}[G]_{r,S})^{r}$ over $\mathbb{Q}[G]$.
Therefore, every element $x\in (\mathbb{Q}\otimes
\bigwedge^{r}_{\mathbb{Z}[G]}U_{S}(F))_{r,S}\cong
\bigwedge^{r}_{\mathbb{Q}[G]_{r,S}}(\mathbb{Q}\otimes
U_{S}(F))_{r,S}$ can be written as $x_{1}\wedge\cdots\wedge x_{r}$,
with $x_{i}\in (\mathbb{Q}\otimes U_{S}(F))_{r,S}$ for all $i$.
\end{rem}
 Let $V=\{v_{1},\cdots, v_{r}\}$ be a set of $r$ places
in $S$, which split completely in $F/K$ and let
$W=\{w_{1},\cdots,w_{r}\}$ be a set of places of $F$ such that
$w_{i}$ lies  above  $v_{i}$, for all $i=1,\cdots, r$. For any place
$w_{i}$, we define the  $G$-equivariant map:
$$
\begin{array}{cccc}
  \lambda_{w_{i}}:& U_{S,\mathcal{T}}(F) & \xymatrix@=2pc{\ar[r]&} & \mathbb{C}[G] \\
 &x & \xymatrix@=2pc{\ar@{|->}[r]&} & -\sum_{\sigma\in
  G}\log(|\sigma(x)|_{w_{i}})\sigma^{-1}.
\end{array}
$$
Rubin's $\mathbb{C}[G]$-linear regulator
$$
\xymatrix@=2pc{\mathrm{Reg}_{S,\mathcal{T}}^{w_{1},\cdots,w_{r}}:
\mathbb{C}\otimes_{\mathbb{Z}}\bigwedge^{r}_{\mathbb{Z}[G]}U_{S,\mathcal{T}}(F)\ar[r]&
\mathbb{C}[G]}
$$
is defined by
$$
\mathrm{Reg}_{S,\mathcal{T}}^{w_{1},\cdots,w_{r}}:=\lambda_{w_{1}}\wedge\cdots\wedge\lambda_{w_{r}}.
$$

Let $\Theta_{S,\mathcal{T}}^{(r)}(0)$ be the coefficient of $s^{r}$
in the Taylor series of $\Theta_{S,\mathcal{T}}$\,;
$$
\Theta_{S,\mathcal{T}}^{(r)}(0):=\displaystyle{\lim_{x\rightarrow
0}} s^{-r}\Theta_{S,\mathcal{T}}^{(r)}(s).
$$
\begin{conj} $\mathbf{RS}(F/K,S,\mathcal{T},r)$.\;Suppose that $(F/K,S,\mathcal{T},r)$ satisfies hypotheses
$\mathbf{H}(F/K,S,\mathcal{T},r)$. Then, for any choice of $V$ and
$W$, there exists a unique element
$\varepsilon_{F,S,\mathcal{T}}\in\Lambda_{S,\mathcal{T}}$
 such that
 $$
\mathrm{Reg}_{S,\mathcal{T}}^{w_{1},\cdots,w_{r}}(\varepsilon_{F,S,\mathcal{T}})=\Theta_{S,\mathcal{T}}^{(r)}(0).
$$
\end{conj}
\begin{rem}
The truth of the Rubin-Stark conjecture
$\mathbf{RS}(F/K,S,\mathcal{T},r)$ does not depend on the particular
choice of $V$ and $W$.
\end{rem}
\vskip7pt Now fix a totally real number field $K$ and let
$r=[K:\mathbb{Q}]$. Recall that
$$
\xymatrix@=2pc{ \chi : G_{K}\ar[r]& \mathcal{O}^{\times}}
$$
is a non-trivial $\overline{\mathbb{Q}}_{p}$-irreducible character,
factoring through a finite abelian  extension $L$ of $K$. Assume
that $L$ and $K_{\infty}$ are linearly disjoint over $K$, and $L$ is
the fixed field of $\ker\chi$. We continue to assume that our
representation is
$$
T=\mathbb{Z}_{p}(1)\otimes \mathcal{O}(\chi^{-1}).
$$
\vskip 6pt
 For a cycle $\mathfrak{r}$ of $K$, let $K(\mathfrak{r})$ be the
maximal $p$-extension inside the ray class field of $K$ modulo
$\mathfrak{r}$. Let $L_{n}=LK_{n}$ and let $\mathfrak{f}_{n}$ denote
the finite part of the conductor of $L_{n}/K$. Remark that
$\mathfrak{f}_{n}$ has the form
$\mathfrak{f}_{n}=\mathfrak{h}\mathfrak{s}_{n}$, where
$\mathfrak{h}$ is prime to $\mathfrak{s}_{n}$ for all $n$ and does
not depend on $n$; moreover $\mathfrak{s}_{n}$ is divisible only by
those prime ideals of $\mathcal{O}_{K}$ which ramify in
$L_{\infty}/L$. For any extension $F$ of $K$, we define
$F(\mathfrak{r})$ as the composite of $K(\mathfrak{r})$ and $F$, and
for any ideal $\mathfrak{a}$, we denote the product of all distinct
prime ideals dividing $\mathfrak{a}$ by $\widetilde{\mathfrak{a}}$.
Let us also assume that $S$ contains the set $\Sigma_{\infty}$ and
at least one finite place, but does not contain any $p$-adic prime
of $K$. For any Galois extension $F$ of $K$, we denote the set of
ramified primes in $F/K$ by $\mathrm{Ram}(F/K)$. Let
$$
S_{F}=S\cup\mathrm{Ram}(F/K).
$$
Following B\"{u}y\"{u}kboduk \cite{Kazim109} we choose
$\mathcal{T}=\{\mathfrak{q}_{0}\}$ where $\mathfrak{q}_{0}$ is a
prime such that $p\nmid \mathbf{N}\mathfrak{q}_{0}-1$ and
$\mathfrak{q}_{0}\nmid 2$, for such
$\mathcal{T}=\{\mathfrak{q}_{0}\}$ we have
$$\widehat{U_{S_{F},\mathcal{T}}(F)}=\widehat{U_{S_{F}}(F)}$$
Moreover, if $F$ is  totally real then the hypothesis
$\mathbf{H}(F/K,S_{F}, \{\mathfrak{q}_{0}\},r)$ is  satisfied. Let
$$
\mathcal{K}_{0}=\{L_{n,\mathfrak{g}},\;
L_{n,\mathfrak{g}}(\mathfrak{r})\;:\; \mathfrak{r} \;\mbox{ is a
finite cycle of $K$ prime to $\mathfrak{q}_{0}\mathfrak{f}_{\chi}p$,
$\mathfrak{g}\mid \widetilde{\mathfrak{h}}$ and $n\in
\mathbb{Z}_{\geq 0}$}\}
$$
where $L_{n,\mathfrak{g}}$ is the maximal subextension of $L_{n}$
whose conductor is prime to
$\widetilde{\mathfrak{h}}\mathfrak{g}^{-1}$, and
$\mathfrak{f}_{\chi}$ is the conductor of $\chi$. In the rest of
this paper we assume \vskip 7pt
 \begin{center}
  \begin{tabbing}
  \hspace{4.5cm}\=\hspace{0.5cm} \= \hspace{0.4cm} \= \kill
  \>$(\mathcal{H}_{0}):$ \> \> \textit{the number field $L$ is totally real,}
\end{tabbing}
\end{center}
and that
\begin{center}
\mbox{\textit{ the conjecture
$\mathbf{RS}(F/K,S_{F},\{\mathfrak{q}_{0}\},r)$ holds, for all $F\in
\mathcal{K}_{0}$.}}
\end{center}
 \begin{rem}
 Let $\varepsilon_{n,\mathfrak{g}}=\varepsilon_{L_{n,\mathfrak{g}},S_{L_{n,\mathfrak{g}}},\{\mathfrak{q}_{0}\}}$
be the Rubin-Stark element for the conjecture\\
$\mathbf{RS}(L_{n,\mathfrak{g}}/K,S_{L_{n,\mathfrak{g}}},\{\mathfrak{q}_{0}\},r)$.
 Since $\mid S_{L_{n,\mathfrak{g}}}\mid> r+1$ for $n\geq 1$, the
element $\varepsilon_{n,\mathfrak{g}}$ lies in
$\mathbb{Q}\otimes_{\mathbb{Z}}\bigwedge^{r}\mathcal{E}_{n,\mathfrak{g}}$,
where $\mathcal{E}_{n,\mathfrak{g}}$ denotes the group of global
units of $L_{n,\mathfrak{g}}$.
 \end{rem}

\begin{deft}\label{Definition of Rubin Strak module}
Let $n$ be a nonnegative integer. We denote by $\mathrm{St}_{n}$ the
$\mathbb{Z}[\mathrm{Gal}(L_{n}/K)]$-module generated by the inverse
images of $\varepsilon_{n,\mathfrak{g}}$ under the map $
\xymatrix@=2pc{ \bigwedge^{r}\mathcal{E}_{n}\ar[r]&
 \mathbb{Q}\otimes\bigwedge^{r}\mathcal{E}_{n}}$
 for all $\mathfrak{g}\mid \widetilde{\mathfrak{h}}$.
\end{deft}\vskip 7pt
Recall that for any number field $F$, Kummer theory gives a
canonical isomorphism
$$
H^{1}(F,\mathbb{Z}_{p}(1))\cong F^{\times,\wedge}:=\varprojlim
F^{\times}/(F^{\times})^{p^{n}}.
$$
Since  $\chi(G_{L_{n}(\mathfrak{r})})=1$ for every $n\geq0$,
$$
H^{1}(L_{n}(\mathfrak{r}),\mathbb{Z}_{p}(1))\otimes
\mathcal{O}(\chi^{-1})\cong
H^{1}(L_{n}(\mathfrak{r}),\mathbb{Z}_{p}(1)\otimes\mathcal{O}(\chi^{-1})).
$$
 Therefore
 \begin{equation}\label{ismophism Lng}
L_{n}(\mathfrak{r})^{\times,\wedge}\otimes
\mathcal{O}(\chi^{-1})\cong
H^{1}(L_{n}(\mathfrak{r}),\mathbb{Z}_{p}(1)\otimes\mathcal{O}(\chi^{-1})).
 \end{equation}
 Let $\varepsilon_{n}(\mathfrak{r})=
 \varepsilon_{L_{n}(\mathfrak{r}),S_{L_{n}(\mathfrak{r})}, \{\mathfrak{q}_{0}\}}$ be the Rubin-Stark
element for
$\mathbf{RS}(L_{n}(\mathfrak{r})/K,S_{L_{n}(\mathfrak{r})},\{\mathfrak{q}_{0}\},r)$.
By Remark \ref{Remark Rubin lattice},
$\varepsilon_{n}(\mathfrak{r})$ can be uniquely written as
$\varepsilon_{1}\wedge\cdots\wedge \varepsilon_{r}$, with
$\varepsilon_{i}\in \mathbb{Q}\otimes L_{n}(\mathfrak{r})^{\times}$.
 Let us note
\begin{equation}\label{definition of chi element of p copletion}
\varepsilon_{n,\chi}(\mathfrak{r}):=\widehat{\varepsilon_{1}}\otimes1_{\chi^{-1}}\wedge\cdots\wedge
\widehat{\varepsilon_{r}}\otimes 1_{\chi^{-1}}
\end{equation}
where $\widehat{\varepsilon_{i}}$ is the image of $\varepsilon_{i}$
by the natural map $ \xymatrix@=1.5pc{ \mathbb{Q}\otimes
L_{n}(\mathfrak{r})^{\times}\ar[r]&
\mathbb{Q}_{p}\otimes_{\mathbb{Z}_{p}} L_{n}(\mathfrak{r})^{\times,
\wedge}}.$
 Then under the isomorphism $(\ref{ismophism Lng})$, we can
view each
$$
\mbox{ $\varepsilon_{n,\chi}(\mathfrak{r})$ as an element of
$\mathbb{Q}_{p}\otimes\bigwedge^{r}H^{1}(L_{n}
(\mathfrak{r}),\mathbb{Z}_{p}(1)\otimes\mathcal{O}(\chi^{-1})).$ }
$$
For any cycle $\mathfrak{r}$ which is prime to
$\mathfrak{q}_{0}\mathfrak{f}_{\chi}p$, and every $n\geq 0$, we
define
\begin{equation}\label{definition of Engr}
c_{n}(\mathfrak{r})=
  \mathrm{cor}_{L_{n+1}(\mathfrak{r}),K_{n}(\mathfrak{r})}^{(r)}
(\varepsilon_{n+1,\chi}(\mathfrak{r})),\quad c_{n}=
  \mathrm{cor}_{L_{n+1},K_{n}}^{(r)}
(\varepsilon_{n+1,\chi})
\end{equation}
where
$\mathrm{cor}_{L_{n+1}(\mathfrak{r}),K_{n}(\mathfrak{r})}^{(r)}$ is
the map
$$
\xymatrix@=2pc{ \mathbb{Q}_{p}\otimes
\bigwedge^{r}H^{1}(L_{n+1}(\mathfrak{r}),T)\ar[r]&
\mathbb{Q}_{p}\otimes \bigwedge^{r}H^{1}(K_{n}(\mathfrak{r}),T)}
$$
induced by the corestrection map
$$
\xymatrix@=2pc{
\mathrm{cor}_{L_{n+1}(\mathfrak{r}),K_{n}(\mathfrak{r})}:
H^{1}(L_{n+1}(\mathfrak{r}),T)\ar[r]& H^{1}(K_{n}(\mathfrak{r}),T)}.
$$
For the convenience of the reader we recall that for any  finite
group $G$ and any  $\mathcal{O}[G]$-module $M$, we have a map
\begin{equation}\label{morphism  iota}
\xymatrix@=2pc{\iota_{M}:
\bigwedge^{r-1}_{\mathcal{O}[G]}\mathrm{Hom}_{\mathcal{O}[G]}(M,\mathcal{O}[G])\ar[r]&
\mathrm{Hom}_{\Phi[G]}(\Phi\otimes\bigwedge^{r}_{\mathcal{O}[G]}M,\Phi\otimes
M)}.
\end{equation}
Indeed, the natural map
$\xymatrix@=2pc{\mathrm{Hom}_{\mathcal{O}[G]}(M,\mathcal{O}[G])\ar[r]&
\mathrm{Hom}_{\Phi[G]}(\Phi\otimes M,\Phi[G])}$ gives a morphism
$$
\xymatrix@=2pc{\bigwedge^{r-1}_{\mathcal{O}[G]}\mathrm{Hom}_{\mathcal{O}[G]}(M,\mathcal{O}[G])\ar[r]&
\bigwedge^{r-1}_{\Phi[G]}\mathrm{Hom}_{\Phi[G]}(\Phi\otimes
M,\Phi[G])}.
$$
On the other hand, the map
$$
\xymatrix@=2pc{f: \mathrm{Hom}_{\Phi[G]}(\Phi\otimes
M,\Phi[G])\ar[r]& \mathrm{Hom}_{\Phi[G]}(\bigwedge^{s}_{\Phi[G]}
\Phi\otimes M,\bigwedge^{s-1}_{\Phi[G]}\Phi\otimes M)}
$$
 defined by
 $f(\psi)(m_{1}\wedge\cdots\wedge
 m_{s})=\displaystyle{\sum_{i=1}^{s}}(-1)^{i+1}\psi(m_{i})m_{1}\wedge\cdots\wedge
m_{i-1}\wedge m_{i+1}\cdots\wedge m_{s},$ and
 iterated from $s=r$ to $s=2$ gives a morphism
\begin{equation*}\label{Morphism de Rubin in power extor}
\xymatrix@=2pc{ \bigwedge^{r-1}_{\Phi[G]}\mathrm{Hom}_{\Phi[G]}(
\Phi\otimes
M,\Phi[G])\ar[r]&\mathrm{Hom}_{\Phi[G]}(\bigwedge^{r}_{\Phi[G]}
\Phi\otimes M,\Phi\otimes M)}.
\end{equation*}
Let us recall also that for any finite Galois extensions $F\subset
F^{\prime}$ of $K$, the norm map from $F^{\prime}$ to $F$ induces a
homomorphism
\begin{equation}\label{map induced by restriction and group}
\xymatrix@=2pc{
\bigwedge^{r-1}_{\mathcal{O}[\Delta_{F^{\prime}}]}\mathrm{Hom}_{\mathcal{O}[\Delta_{F^{\prime}}]}(H^{1}(F^{\prime},T),\mathcal{O}[\Delta_{F^{\prime}}])\ar[r]
&
\bigwedge^{r-1}_{\mathcal{O}[\Delta_{F}]}\mathrm{Hom}_{\mathcal{O}[\Delta_{F}]}(H^{1}(F,T),\mathcal{O}[\Delta_{F}])}
\end{equation}
where $\Delta_{F}=\mathrm{Gal}(F/K)$. In particular, the collection
$$\big(\bigwedge^{r-1}\mathrm{Hom}_{\mathcal{O}[\Delta_{K_{n}(\mathfrak{r})}]}
(H^{1}(K_{n}(\mathfrak{r}),T),\mathcal{O}[\Delta_{K_{n}(\mathfrak{r})}])\big)_{n,\mathfrak{r}}$$
is a projective system for the maps given in $(\ref{map induced by
restriction and group})$.
\begin{deft}\label{definition of psi euler system}
Let $\Psi=\{\psi_{n,\mathfrak{r}}\}_{n,\mathfrak{r}}$ be an
arbitrary element of
$$\displaystyle{\varprojlim_{n,\mathfrak{r}}}\bigwedge^{r-1}\mathrm{Hom}_{\mathcal{O}[\Delta_{K_{n}(\mathfrak{r})}]}
(H^{1}(K_{n}(\mathfrak{r}),T),\mathcal{O}[\Delta_{K_{n}(\mathfrak{r})}]).$$
Let us identify $\psi_{n,\mathfrak{r}}$ with its image
$\iota_{M}(\psi_{n,\mathfrak{r}})$ under $(\ref{morphism iota})$,
$M=H^{1}(K_{n}(\mathfrak{r}),T)$. We define
$$
\varepsilon_{n,\Psi}(\mathfrak{r})^{\chi}:=\psi_{n,\mathfrak{r}}(c_{n}(\mathfrak{r}))\quad\mbox{and}\quad
\varepsilon_{n,\Psi}^{\chi}:=\psi_{n}(c_{n}),
$$
where $c_{n}(\mathfrak{r})$ is given  in Definition
$(\ref{definition of Engr})$.
\end{deft}
By the defining
 integrality property of the elements
 $\varepsilon_{n}(\mathfrak{r})$ and
 Corollary $1.3$ in \cite{Rubin96}, we get
 $$
\varepsilon_{n,\Psi}(\mathfrak{r})^{\chi}\in
H^{1}(K_{n}(\mathfrak{r}),T)\quad\mbox{and}\quad
\varepsilon_{n,\Psi}^{\chi}\in H^{1}(K_{n},T) .
$$
Remark that for every $n\geq 0$, we have
$\mathrm{cor}_{K_{n}(1),K_{n}}(\varepsilon_{n,\Psi}(\mathfrak{1})^{\chi})=\varepsilon_{n,\Psi}^{\chi}$.
\begin{pro}\label{Euler system}
The collection $\{\varepsilon_{n,\Psi}(\mathfrak{r})^{\chi}\}_{n\geq
0, \mathfrak{r}}$ is an Euler system for the
$\mathrm{G}_{K}$-representation $T$, in the sense of
\cite[Definition 2.1.1]{Rubin00}.
\end{pro}
\noindent \textbf{Proof.} This is a consequence of Proposition $6.2$
in \cite{Rubin96}.\hfill $\square$
\section{\bf Modifying the local condition at $p$.}
In this section, we modify the classical local conditions at the
primes above $p$ to obtain a Selmer structure $\mathcal{L}$ on $T$.
To this end we assume  the hypotheses $(\mathcal{H}_{1})$ and
$(\mathcal{H}_{2})$ all throughout. The following theorem is crucial
for our purpose. It is a direct consequence of \cite[Theorem
2.2]{Greither96}.
\begin{theo}\label{Pro H1IW rank }
For any $p$-adic prime $v$ of $K(\mathfrak{r})$, the
$\mathcal{O}[[\mathrm{Gal}(K_{\infty}(\mathfrak{r})_{v}/\mathbb{Q}_{p})]]$-module
$$
H^{1}_{Iw}(K(\mathfrak{r})_{v},T):=\varprojlim_{n}
H^{1}(K_{n}(\mathfrak{r})_{v},T)
$$
is free of rank one.
\end{theo}
\noindent \textbf{Proof.} Let $F=L(\mathfrak{r})$ and
$F_{n}=L_{n}(\mathfrak{r})$. Fix a place $w$ of $F$ lying above $v$.
By \cite[Theorem 2.2]{Greither96}, we have an exact sequence of
$\mathbb{Z}_{p}[[\mathrm{Gal}(F_{w}(\mu_{p^{\infty}})/\mathbb{Q}_{p})]]$-modules
$$
\xymatrix@=2pc{ 0\ar[r]& \mathbb{Z}_{p}(1)\ar[r]& \varprojlim_{n}
U^{1}(F_{n,w}(\mu_{p}))\ar[r]^-{h}& \mathcal{V}(1)\ar[r]&
\mathbb{Z}_{p}(1)\ar[r]&0}
$$
where the module $\mathcal{V}(1)$ is free cyclic over
$\mathbb{Z}_{p}[[\mathrm{Gal}(F_{w}(\mu_{p^{\infty}})/\mathbb{Q}_{p})]]$.
Taking the
$\mathrm{Gal}(F_{w}(\mu_{p^{\infty}})/F_{\infty,w})$-cohomology of
the exact sequences
$$
\xymatrix@=2pc{ 0\ar[r]& \mathbb{Z}_{p}(1)\ar[r]& \varprojlim_{n}
U^{1}(F_{n,w}(\mu_{p}))\ar[r]^-{h}& \mathrm{im}(h)\ar[r]& 0},
$$
$$
\xymatrix@=2pc{ 0\ar[r]& \mathrm{im}(h)\ar[r]& \mathcal{V}(1)\ar[r]&
\mathbb{Z}_{p}(1)\ar[r]&0}
$$
and remarking that $$ H^{i}(F_{w}(\mu_{p^{\infty}})/F_{\infty,w},
\mathbb{Z}_{p}(1))=0, \quad \mbox{for $i\geq 0$}
$$
we obtain
$$
H^{0}(F_{w}(\mu_{p^{\infty}})/F_{\infty,w},
\varprojlim_{n}U^{1}(F_{n,w}(\mu_{p}))\cong
H^{0}(F_{w}(\mu_{p^{\infty}})/F_{\infty,w}, \mathcal{V}(1))
$$
It follows that
$$
 \varprojlim_{n}U^{1}(F_{n,w})\cong
H^{0}(F_{w}(\mu_{p^{\infty}})/F_{\infty,w}, \mathcal{V}(1)).
$$
Since $\chi(\mathrm{Frob}_{v})\neq 1$, we obtain

\begin{eqnarray*}
  H^{1}_{Iw}(K(\mathfrak{r})_{v}, T) &\cong &
(\varprojlim_{n}U^{1}(F_{n,w})\otimes\mathcal{O}(\chi^{-1}))^{\mathrm{Gal}(F_{\infty,w}/K_{\infty}(\mathfrak{r})_{v})}  \\
   &\cong& (H^{0}(F_{w}(\mu_{p^{\infty}})/F_{\infty,w},
\mathcal{V}(1))\otimes\mathcal{O}(\chi^{-1}))^{\mathrm{Gal}(F_{\infty,w}/K_{\infty}(\mathfrak{r})_{v})}.
\end{eqnarray*}
Then $ H^{1}_{Iw}(K(\mathfrak{r})_{v}, T)$ is a free
$\mathcal{O}[[\mathrm{Gal}(K_{\infty}(\mathfrak{r})_{v}/\mathbb{Q}_{p})]]$-module
of rank one. \hfill $\square$
\begin{coro}\label{Lemma structur of H(Kp)}
The $\mathcal{O}[[\mathrm{Gal}(K_{\infty}(\mathfrak{r})/K)]]$-module
$$
H^{1}_{Iw}(K(\mathfrak{r})_{p},T)
$$
is free of rank $[K:\mathbb{Q}]$. In particular,
$H^{1}_{Iw}(K_{p},T)$ is a free $\Lambda$-module of rank
$[K:\mathbb{Q}]$.
\end{coro}
\noindent\textbf{Proof.} Let $v$ be a $p$-adic prime of
$K(\mathfrak{r})$. By  Theorem \ref{Pro H1IW rank } we see that
$H^{1}_{Iw}(K(\mathfrak{r})_{v},T)$ is a free
$\mathcal{O}[[\mathrm{Gal}(K_{\infty}(\mathfrak{r})_{v}/K_{w})]]$-module
of rank $[K_{w}:\mathbb{Q}_{p}]$, where $w$ is a place of $K$ lying
below $v$. Then
$$
H^{1}_{Iw}(K(\mathfrak{r})_{p},T)=\oplus_{v\mid
p}H^{1}_{Iw}(K(\mathfrak{r})_{v},T)
$$
 is a free
$\mathcal{O}[[\mathrm{Gal}(K_{\infty}(\mathfrak{r})/K)]]$-module of
rank $\displaystyle{\sum_{w\mid
p}}[K_{w}:\mathbb{Q}_{p}]=[K:\mathbb{Q}] $.\hfill $\square$\vskip
7pt Let $\mathcal{K}=\bigcup_{n,\mathfrak{r}}K_{n}(\mathfrak{r})$
and let $\mathcal{G}$ denote the Galois group
$\mathrm{Gal}(\mathcal{K}/K)$.
\begin{coro}\label{coro the rank of V as Galois module}
The $\mathcal{O}[[\mathcal{G}]]$-module
$\mathbb{V}:=\varprojlim_{n,\mathfrak{r}}H^{1}(K_{n}(\mathfrak{r})_{p},T)=
\varprojlim_{\mathfrak{r}}H^{1}_{Iw}(K(\mathfrak{r})_{p},T)$ is free
of rank $[K:\mathbb{Q}]$.
\end{coro}
\noindent\textbf{Proof.} Immediate after Corollary \ref{Lemma
structur of H(Kp)}. \hfill $\square$\vskip 6pt
\begin{rem}
This corollary is a generalization of \cite[Corollary
3.10]{Kazim109}.
\end{rem}
If $F\subset F^{\prime}$ is an extension of fields, we will write
$F\subset_{f}F^{\prime}$ to indicate that $[F^{\prime}:F]$ is
finite.
\begin{pro}\label{pro  coinvariant and cokernel}
Let $K\subset_{f}F\subset \mathcal{K}$ be a finite extension and let
$\mathcal{G}_{F}$ denote the Galois group
$\mathrm{Gal}(\mathcal{K}/F)$. Let
$d=\displaystyle{\max_{\mathfrak{p}\mid
p}}\{v_{p}(1-\chi(\mathrm{Frob}_{\mathfrak{p}}))\}$, where
$\mathfrak{p}$ is a prime of $K$. Then the canonical map
$$
\xymatrix@=2pc{ \mathbb{V}_{\mathcal{G}_{F}}\ar[r]& H^{1}(F_{p},T)}
$$
is injective, with cokernel annihilated by $p^{d}$.

\end{pro}
 \noindent \textbf{Proof.} For any $p$-adic place $v$ of $F$, we
fixe a place $\overline{v}$ of $\mathcal{K}$ lying above $v$.
%
One has
$$
\mathbb{V}_{\mathcal{G}_{F}}\cong\bigoplus_{\substack {v\mid p\\
v\subset F}} (\varprojlim_{F_{v}\subset_{f}F^{\prime}_{w}\subset
\mathcal{K}_{\overline{v}}}H^{1}(F^{\prime}_{w},T))_{D_{v}(\mathcal{K}/F)}
$$
where $D_{v}(\mathcal{K}/F)$ is the decomposition subgroup of $v$ in
$\mathcal{K}/F$. By dualizing the inflation-restriction exact
sequence
$$
\xymatrix@=2pc{H^{1}(D_{v}(\mathcal{K}/F),(T^{\ast})^{G_{\mathcal{K}_{\overline{v}}}})\ar@{^{(}->}[r]&H^{1}(F_{v},T^{\ast})\ar[r]&
H^{1}(\mathcal{K}_{\overline{v}},T^{\ast})^{D_{v}(\mathcal{K}/F)}\ar[r]&H^{2}(D_{v}(\mathcal{K}/F),(T^{\ast})^{G_{\mathcal{K}_{\overline{v}}}})},
$$
we obtain an exact sequence
$$
\xymatrix@=1pc{H^{2}(D_{v}(\mathcal{K}/F),(T^{\ast})^{G_{\mathcal{K}_{\overline{v}}}})^{\vee}\ar[r]&
(H^{1}(\mathcal{K}_{\overline{v}},T^{\ast})^{D_{v}(\mathcal{K}/F)})^{\vee}\ar[r]&
H^{1}(F_{v},T^{\ast})^{\vee}\ar@{->>}[r]&
(H^{1}(D_{v}(\mathcal{K}/F),(T^{\ast})^{G_{\mathcal{K}_{\overline{v}}}}))^{\vee}}
$$
Since
$H^{2}(D_{v}(\mathcal{K}/F),(T^{\ast})^{G_{\mathcal{K}_{\overline{v}}}})^{\vee}$
is a torsion $\mathcal{O}$-module and
$$
(H^{1}(\mathcal{K}_{\overline{v}},T^{\ast})^{D_{v}(\mathcal{K}/F)})^{\vee}\cong
(\displaystyle{\varprojlim_{ F_{v}\subset_{f}F^{\prime}_{w}\subset
\mathcal{K}_{\overline{v}}}}H^{1}(F^{\prime}_{w},T))_{D_{v}(\mathcal{K}/F)}
$$
is a torsion-free $\mathcal{O}$-module, then we have an exact
sequence:
$$
\xymatrix@=2pc{0\ar[r]&
(H^{1}(\mathcal{K}_{\overline{v}},T^{\ast})^{D_{v}(\mathcal{K}/F)})^{\vee}\ar[r]&
H^{1}(F_{v},T^{\ast})^{\vee}\ar@{->>}[r]&
(H^{1}(D_{v}(\mathcal{K}/F),(T^{\ast})^{G_{\mathcal{K}_{\overline{v}}}}))^{\vee}}.
$$
Since $T^{\ast}=D(\chi)$, then for every  $ K\subset F\subset
\mathcal{K}$, $ H^{0}(K_{p},T^{\ast})=H^{0}(F_{p},T^{\ast}).$
Therefore
$$
p^{d}.(T^{\ast})^{G_{\mathcal{K}_{\overline{v}}}}=0
$$
where $d=\displaystyle{\max_{\mathfrak{p}\mid
p}}\{v_{p}(1-\chi(\mathrm{Frob}_{\mathfrak{p}}))\}$, $\mathfrak{p}$
is a prime of $K$. It follows that the canonical map
$$
\xymatrix@=2pc{ \mathbb{V}_{\mathcal{G}_{F}}\ar[r]& H^{1}(F_{p},T)}
$$
is injective, with cokernel annihilated by $p^{d}$. \hfill
$\square$\vskip 6pt As  B\"{u}y\"{u}kboduk did in \cite{Kazim109},
we fix a free $\mathcal{O}[[\mathrm{Gal}(\mathcal{K}/K)]]$-direct
summand $\mathbb{L}$ inside of
$$\mathbb{V}=\displaystyle{\varprojlim_{K\subset_{f}F\subset
\mathcal{K}}}H^{1}(F_{p},T)= \varprojlim_{K\subset_{f}F\subset
\mathcal{K}}\mathbb{V}_{\mathcal{G}_{F}}$$ which is free of rank one
as $\mathcal{O}[[\mathrm{Gal}(\mathcal{K}/K)]]$-module.
 Recall that $\Sigma$ is a finite set of
places of $K$ containing all infinite places, all $p$-adic places
and all places where $T$ is ramified.
\begin{deft}\label{Definition of modified Selmer}
 Define the modified Selmer
structure $\mathcal{L}$ on $T$ by
\begin{itemize}
    \item $\Sigma(\mathcal{L})=\Sigma$,
    \item if $w\nmid p$,
    $H^{1}_{\mathcal{L}}(F_{w},T)=H^{1}_{\mathcal{F}_{can}}(F_{w},T)$,
    \item   $H^{1}_{\mathcal{L}}(F_{p},T)\subset H^{1}(F_{p},T)$ as the
$\mathcal{O}$-saturation of \,$\mathbb{L}_{\mathcal{G}_{F}}$ in
$H^{1}(F_{p},T)$.
\end{itemize}
\end{deft}
\subsection{\bf Choosing homomorphisms} Let us keep the same notation as above. In this subsection, we show
the existence of a homomorphism
$\Psi^{\prime}=(\Psi^{\prime}_{F})_{K\subset_{f}F\subset
\mathcal{K}}$ such that
$$\Psi^{\prime}_{F}\big(\bigwedge^{r}H^{1}(F_{p},T)\big)\subset
H^{1}_{\mathcal{L}}(F_{p},T).$$ Fix a basis
$\{\psi_{\mathbb{L}}^{(i)}\}_{i=1}^{r-1}$ of the free
$\mathcal{O}[[\mathrm{Gal}(\mathcal{K}/K)]]$-module $$
\mathrm{Hom}_{\mathcal{O}[[\mathrm{Gal}(\mathcal{K}/K)]]}(\mathbb{V}/\mathbb{L},\mathcal{O}[[\mathrm{Gal}(\mathcal{K}/K)]])
$$ of rank $r-1$. This in return fixes a basis
$\{\psi_{\mathcal{L}_{F}}^{(i)}\}_{i=1}^{r-1}$ for the free
$\mathcal{O}[\Delta_{F}]$-module
$$
\mathrm{Hom}_{\mathcal{O}[\Delta_{F}]}(\mathbb{V}_{\mathcal{G}_{F}}/\mathbb{L}_{\mathcal{G}_{F}},
\mathcal{O}[\Delta_{F}])
$$
for all $K\subset_{f}F\subset \mathcal{K}$, such that the
homomorphisms $\{\psi_{\mathcal{L}_{F}}^{(i)}\}_{i=1}^{r-1}$ are
compatible with respect to the maps
$$
\xymatrix@=2pc{\mathrm{Hom}_{\mathcal{O}[\Delta_{F^{\prime}}]}
(\mathbb{V}_{\mathcal{G}_{F^{\prime}}}/\mathbb{L}_{\mathcal{G}_{F^{\prime}}},
\mathcal{O}[\Delta_{F^{\prime}}])\ar[r]&
\mathrm{Hom}_{\mathcal{O}[\Delta_{F}]}(\mathbb{V}_{\mathcal{G}_{F}}/\mathbb{L}_{\mathcal{G}_{F}},
\mathcal{O}[\Delta_{F}])}
$$
for $F\subset_{f}F^{\prime}$. Let $\psi_{F}^{(i)}$ denote the image
of $\psi_{\mathcal{L}_{F}}^{(i)}$ under the canonical injection
$$
\xymatrix@=2pc{
\mathrm{Hom}_{\mathcal{O}[\Delta_{F}]}(\mathbb{V}_{\mathcal{G}_{F}}/\mathbb{L}_{\mathcal{G}_{F}},
\mathcal{O}[\Delta_{F}])\ar@{^{(}->}[r]&\mathrm{Hom}_{\mathcal{O}[\Delta_{F}]}(\mathbb{V}_{\mathcal{G}_{F}},
\mathcal{O}[\Delta_{F}])}.
$$
Remark that the map
$$
\xymatrix@=2pc{\Psi_{\mathcal{L}_{F}}:=\bigoplus_{i=1}^{r-1}\psi_{F}^{(i)}:
\mathbb{V}_{\mathcal{G}_{F}}\ar[r]& \mathcal{O}[\Delta_{F}]^{r-1}}
$$
is surjective and
$\ker(\Psi_{\mathcal{L}_{F}})=\mathbb{L}_{\mathcal{G}_{F}}$. \vskip
6pt Define
$$
\Psi_{F}:=\psi_{F}^{(1)}\wedge \psi_{F}^{(2)}\wedge\cdots\wedge
\psi_{F}^{(r-1)}\in
\bigwedge^{r-1}\mathrm{Hom}_{\mathcal{O}[\Delta_{F}]}(\mathbb{V}_{\mathcal{G}_{F}},\mathcal{O}[\Delta_{F}]).
$$
We may therefore regard $\Psi:=(\Psi_{F})_{K\subset_{f}F\subset
\mathcal{K}}$ as an element of the module
$$
\varprojlim_{K\subset_{f}F\subset
\mathcal{K}}\bigwedge^{r-1}\mathrm{Hom}_{\mathcal{O}[\Delta_{F}]}(\mathbb{V}_{\mathcal{G}_{F}},\mathcal{O}[\Delta_{F}])
$$
\begin{pro}\label{pro Psi and coinvariant}
Let $\Psi:=(\Psi_{F})_{K\subset_{f}F\subset \mathcal{K}}$ be as
above. Then for every $K\subset_{f}F\subset \mathcal{K}$, $\Psi_{F}$
induces an isomorphism
$$
\xymatrix@=2pc{ \Psi_{F}:
\bigwedge^{r}\mathbb{V}_{\mathcal{G}_{F}}\ar[r]^-{\sim}&
\ker(\Psi_{\mathcal{L}_{F}})=\mathbb{L}_{\mathcal{G}_{F}}}
$$
\end{pro}
\noindent \textbf{Proof.} The proof is identical to the proof of
\cite[Proposition 3.17]{Kazim109}, which follows the proof of Lemma
$\mathrm{B}.1$ of  \cite{MR04}.\hfill $\square$
\begin{pro}\label{pro element psi0}
There exists an element
$$\Psi^{\prime}=(\psi^{\prime}_{F})_{K\subset_{f}F\subset \mathcal{K}}\in
\varprojlim_{K\subset_{f}F\subset
\mathcal{K}}\bigwedge^{r-1}\mathrm{Hom}_{\mathcal{O}[\Delta_{F}]}(H^{1}(F_{p},T),\mathcal{O}[\Delta_{F}])$$
such that for any $K\subset_{f}F\subset \mathcal{K}$,
$$
\psi^{\prime}_{F}(\bigwedge^{r}H^{1}(F_{p},T))\subset
H^{1}_{\mathcal{L}}(F_{p},T).
$$
\end{pro}
\noindent \textbf{Proof.} Thanks to Proposition $\ref{pro
coinvariant and cokernel}$, the canonical map
$$
\xymatrix@=2pc{ \mathbb{V}_{\mathcal{G}_{F}}\ar[r]& H^{1}(F_{p},T)}
$$
is injective, with cokernel annihilated by $p^{d}$. Then the
cokernel of the canonical map
$$
\xymatrix@=2pc{
\mathrm{Hom}_{\mathcal{O}[\Delta_{F}]}(H^{1}(F_{p},T),\mathcal{O}[\Delta_{F}])\ar[r]&
\mathrm{Hom}_{\mathcal{O}[\Delta_{F}]}(\mathbb{V}_{\mathcal{G}_{F}},\mathcal{O}[\Delta_{F}])}
$$
is annihilated by $p^{d}$. Hence the cokernel of
\begin{equation}\label{map H1 vers coinvariant}
\xymatrix@=2pc{ \displaystyle{\varprojlim_{K\subset_{f}F\subset
\mathcal{K}}}\bigwedge^{r-1}\mathrm{Hom}_{\mathcal{O}[\Delta_{F}]}(H^{1}(F_{p},T),\mathcal{O}[\Delta_{F}])\ar[r]&
\displaystyle{\varprojlim_{K\subset_{f}F\subset
\mathcal{K}}}\bigwedge^{r-1}
\mathrm{Hom}_{\mathcal{O}[\Delta_{F}]}(\mathbb{V}_{\mathcal{G}_{F}},\mathcal{O}[\Delta_{F}])}
\end{equation}
is annihilated by $p^{d(r-1)}$. Therefore $p^{d(r-1)}\Psi$ is an
element of the image of the map $(\ref{map H1 vers coinvariant})$,
where $\Psi$ is defined in Proposition $\ref{pro Psi and
coinvariant}$. Then
$$
(p^{d(r-1)}\psi_{F})(\bigwedge^{r}H^{1}(F_{p},T))\subset
\mathbb{L}_{\mathcal{G}_{F}}\subset
H^{1}_{\mathcal{L}}(F_{p},T).\quad \square
$$
\vskip 7pt
\subsection{\bf Kolyvagin systems for $(T,\mathcal{L})$} In this
subsection, we show that the Kolyvagin's derivative class associated
to the Euler system of Rubin-Stark elements defines a Kolyvagin
system for the modified Selmer structure $\mathcal{L}$.\vskip 6pt
 Let $c_{K,\infty}=\{\varepsilon_{n,\Psi}^{\chi}\}\in
 H^{1}_{Iw}(K,T)$ denote the element corresponding to the Euler
 system
 $\{\varepsilon_{n,\Psi}(\mathfrak{r})^{\chi}\}_{n,\mathfrak{r}}$ in
 $H^{1}_{Iw}(K,T)=\displaystyle{\varprojlim_{n}}H^{1}(K_{n},T)$.
\begin{pro}\label{Proposition of locp}
Let $\mathbb{H}_{\mathbb{L}}$ be the set of the maps
$\Psi=(\Psi_{F})_{F}\in \displaystyle{\varprojlim_{K\subset F\subset
\mathcal{K}}}\bigwedge^{r-1}\mathrm{Hom}_{\mathcal{O}[\Delta_{F}]}(\mathbb{V}_{\mathcal{G}_{F}},
\mathcal{O}[\Delta_{F}])$ such that $
\Psi_{F}(\bigwedge^{r}\mathbb{V}_{\mathcal{G}_{F}})=\mathbb{L}_{\mathcal{G}_{F}}$.
Let $\mathrm{loc}_{p}$ denote the localization map at $p$;
$$
\xymatrix@=2pc{\mathrm{loc}_{p}\;: H^{1}(K,T)\ar[r]&H^{1}(K_{p},T)}.
$$
Then
$$
\{\mathrm{loc}_{p}(\varepsilon_{0,\Psi}^{\chi}): \Psi\in
\mathbb{H}_{\mathbb{L}} \}=
[\bigwedge^{r}\mathbb{V}_{\mathcal{G}}:\mathcal{O}.\mathrm{loc}^{(r)}_{p}(c_{0})]\mathbb{L}_{\mathcal{G}}.
$$
\end{pro}
\noindent \textbf{Proof.} The proof is identical to the proof of
Corollary $3.5$ of \cite{Kazim108} line by line. \hfill
$\square$\vskip 6pt
\begin{rem}\label{remark indice fini} If the localization map
$\xymatrix@=2pc{\mathrm{loc}_{p} : H^{1}(K,T)\ar[r]&
H^{1}(K_{p},T)}$ is injective, then by \cite[Proposition 6.6
(ii)]{Rubin96}, we see that
$[\bigwedge^{r}\mathbb{V}_{\mathcal{G}}:\mathcal{O}.\mathrm{loc}^{(r)}_{p}(c_{0})]<\infty$.
\end{rem}
 Let $F$ be a finite extension of $K$ in $\mathcal{K}$ and
let $\mathcal{F}$ be a Selmer structure on $T$. For any cycle
$\mathfrak{r}$ of $K$, we write $\mathcal{F}^{\mathfrak{r}}$ for the
Selmer structure defined by
\begin{itemize}
    \item
$\Sigma(\mathcal{F}^{\mathfrak{r}})=\Sigma(\mathcal{F})\cup
\Sigma_{\mathfrak{r}}$
    \item $H^{1}_{\mathcal{F}^{\mathfrak{r}}}(F_{w},T)= \left\{
            \begin{array}{ll}
            H^{1}_{\mathcal{F}}(F_{w},T)  , & \hbox{if $w\in \Sigma(\mathcal{F})-\Sigma_{\mathfrak{r}}$;} \\
             H^{1}(F_{w},T) , & \hbox{$w\in\Sigma_{\mathfrak{r}}$.}
            \end{array}
          \right.
$
\end{itemize}
where $\Sigma_{\mathfrak{r}}=\{w\subset F\;; w\mid
\mathfrak{r}\}.$\\
Let $M$ be a power of a uniformizer of $\mathcal{O}$ and let
$W_{M}=T/MT$. Recall that for any Euler system $\mathbf{c}$ of $T$
we can associate a Kolyvagin derivative class
$\kappa_{[K_{n},\mathfrak{r},M]}$, see \cite[\S 4.4]{Rubin00}.
Recall also that
$$
\kappa_{[K_{n},\mathfrak{r},M]}\in
H^{1}_{\mathcal{F}_{can}^{\mathfrak{r}}}(K_{n},W_{M})
$$
 c.f.\,\cite[Theorem 4.5.1]{Rubin00}. Next, we construct an Euler
system $\mathbf{c}$ of $T$ such that
$$
\kappa_{[K_{n},\mathfrak{r},M]}\in
H^{1}_{\mathcal{L}^{\mathfrak{r}}}(K_{n},W_{M}).
$$
For this we need some  facts about the local  condition
$\mathcal{L}$.
\begin{lem}\label{pd s local condition}

 Let
$\mathcal{L}_{F}:=\mathbb{L}_{\mathcal{G}_{F}}$, then
$$
p^{d}.(H^{1}_{\mathcal{L}}(F_{p},T)/\mathcal{L}_{F})=0
$$
where $d=\displaystyle{\max_{\mathfrak{p}\mid
p}}\{v_{p}(1-\chi(\mathrm{Frob}_{\mathfrak{p}}))\}$.
\end{lem}
\noindent \textbf{Proof.} First consider, the exact sequence
$$
\xymatrix@=2pc{0\ar[r]& H^{1}_{\mathcal{L}}(F_{p},T)/\mathcal{L}_{F}
\ar[r]& H^{1}(F_{p},T)/\mathcal{L}_{F}\ar@{->>}[r]&
H^{1}(F_{p},T)/H^{1}_{\mathcal{L}}(F_{p},T)}
$$
By definition, the $\mathcal{O}$-module
$H^{1}(F_{p},T)/H^{1}_{\mathcal{L}}(F_{p},T)$ is torsion-free, then
$$
\mathrm{tor}_{\mathcal{O}}(H^{1}(F_{p},T)/\mathcal{L}_{F})=
H^{1}_{\mathcal{L}}(F_{p},T)/\mathcal{L}_{F}.
$$
Second the facts that  $\mathbb{V}_{\mathcal{G}_{F}}
/\mathcal{L}_{F}$ is $\mathcal{O}$-torsion-free and
$H^{1}(F_{p},T)/\mathbb{V}_{\mathcal{G}_{F}}$ is
$\mathcal{O}$-torsion, and the exact sequence
$$
\xymatrix@=2pc{0\ar[r]& \mathbb{V}_{\mathcal{G}_{F}}
/\mathcal{L}_{F}\ar[r]& H^{1}(F_{p},T)/\mathcal{L}_{F}\ar@{->>}[r]&
H^{1}(F_{p},T)/\mathbb{V}_{\mathcal{G}_{F}}}
$$
show that
$$
 \mathrm{tor}_{\mathcal{O}}(H^{1}(F_{p},T)/\mathcal{L}_{F})\cong
 H^{1}(F_{p},T)/\mathbb{V}_{\mathcal{G}_{F}}.
$$
 Then, by Proposition \ref{pro coinvariant and cokernel}, we get
$$
p^{d}.(H^{1}_{\mathcal{L}}(F_{p},T)/\mathcal{L}_{F})=0. \quad
\square
$$
\begin{pro}\label{Pro Cokernel of WM and power of p}
 Let
$G_{n,\mathfrak{r}}=\mathrm{Gal}(K_{n}(\mathfrak{r})/K_{n})$. Then
the cokernel of
$$
\xymatrix@=2pc{ H^{1}_{\mathcal{L}}(K_{n,p},W_{M})\ar[r]&
H^{1}_{\mathcal{L}}(K_{n}(\mathfrak{r})_{p},W_{M})^{G_{n,\mathfrak{r}}}}
$$
is annihilated by $p^{d}$.
\end{pro}
\noindent \textbf{Proof.} Let $F=K_{n}(\mathfrak{r})$ or $F=K_{n}$.
Consider the exact sequence
$$
\xymatrix@=2pc{0\ar[r]&T\ar[r]^-{M}&T\ar[r]&W_{M}\ar[r]&0}.
$$
 Using \cite[Lemma 3.7.1]{MR04}, \cite[Lemma 1.1.5]{MR04} and the
fact that $H^{1} (F_{p},T)$ is torsion-free $\mathcal{O}$-module, we
get an exact sequence
$$
\xymatrix@=2pc{0\ar[r]&H^{1}_{\mathcal{L}}(F_{p},T)\ar[r]^-{M}&H^{1}_{\mathcal{L}}(F_{p},T)\ar[r]&H^{1}_{\mathcal{L}}(F_{p},W_{M})\ar[r]&0}.
$$
Since the restriction $\xymatrix@=1.5pc{\mathrm{res}:
H^{1}(K_{n,p},T)\ar[r]&
H^{1}(K_{n}(\mathfrak{r})_{p},T)^{G_{n,\mathfrak{r}}}}$ is an
isomorphism and the $\mathcal{O}$-module
$H^{1}(K_{n,p},T)/H^{1}_{\mathcal{L}}(K_{n,p},T)$ is torsion-free,
the commutative diagram
$$
\xymatrix@=1.5 pc{ 0\ar[r]&
H^{1}_{\mathcal{L}}(K_{n,p},T)\ar[r]\ar[d]^-{res}&
H^{1}(K_{n,p},T)\ar[r]\ar[d]^-{res}_-{\wr}&
H^{1}(K_{n,p},T)/H^{1}_{\mathcal{L}}(K_{n,p},T)\ar[r]\ar[d]& 0\\
0\ar[r]&
H^{1}_{\mathcal{L}}(K_{n}(\mathfrak{r})_{p},T)^{G_{n,\mathfrak{r}}}\ar[r]&
H^{1}(K_{n}(\mathfrak{r})_{p},T)^{G_{n,\mathfrak{r}}}\ar[r]&
(H^{1}(K_{n}(\mathfrak{r})_{p},T)/H^{1}_{\mathcal{L}}(K_{n}(\mathfrak{r})_{p},T))^{G_{n,\mathfrak{r}}}&
}
$$
shows that the map  $\xymatrix@=2pc{
H^{1}_{\mathcal{L}}(K_{n,p},T)\ar[r]^-{\mathrm{res}}&
H^{1}_{\mathcal{L}}(K_{n}(\mathfrak{r})_{p},T)^{G_{n,\mathfrak{r}}}}$
is an isomorphism. Therefore, we have an exact commutative diagram
$$
\xymatrix@=1.5pc{0\ar[r]&H^{1}_{\mathcal{L}}(K_{n,p},T)\ar[r]^-{M}\ar[d]^-{res}_-{\wr}&
H^{1}_{\mathcal{L}}(K_{n,p},T)\ar[r]\ar[d]^-{res}_-{\wr}&
H^{1}_{\mathcal{L}}(K_{n,p},W_{M})\ar[d]\ar[r]&0\\
0\ar[r]&H^{1}_{\mathcal{L}}(K_{n}(\mathfrak{r})_{p},T)^{G_{n,\mathfrak{r}}}\ar[r]^-{M}&
H^{1}_{\mathcal{L}}(K_{n}(\mathfrak{r})_{p},T)^{G_{n,\mathfrak{r}}}\ar[r]&
H^{1}_{\mathcal{L}}(K_{n}(\mathfrak{r})_{p},W_{M})^{G_{n,\mathfrak{r}}}&
}
$$
Then by the snake lemma, we get
\begin{eqnarray*}
  \mathrm{coker}(\xymatrix@=1pc{ H^{1}_{\mathcal{L}}(K_{n,p},W_{M})\ar[r]&
   H^{1}_{\mathcal{L}}(K_{n}(\mathfrak{r})_{p},W_{M})^{G_{n,\mathfrak{r}}}}) &\cong&
   \mathrm{coker}(\xymatrix@=1pc{ H^{1}_{\mathcal{L}}(K_{n}(\mathfrak{r})_{p},T)^{G_{n,\mathfrak{r}}}
   \ar[r]& H^{1}_{\mathcal{L}}(K_{n}(\mathfrak{r})_{p},W_{M})^{G_{n,\mathfrak{r}}}})\\
   &=& H^{1}(G_{n,\mathfrak{r}},H^{1}_{\mathcal{L}}(K_{n}(\mathfrak{r})_{p},T))[M]
\end{eqnarray*}
where
$H^{1}(G_{n,\mathfrak{r}},H^{1}_{\mathcal{F}}(K_{n}(\mathfrak{r})_{p},T))[M]$
is the submodule of
$H^{1}(G_{n,\mathfrak{r}},H^{1}_{\mathcal{F}}(K_{n}(\mathfrak{r})_{p},T))$
annihilated by $M$. Therefore  it suffices to prove that
$$
p^{d}.H^{1}(G_{n,\mathfrak{r}},H^
{1}_{\mathcal{F}}(K_{n}(\mathfrak{r})_{p},T))=0.
$$
For this, consider the exact sequence
$$
\xymatrix@=2pc{0\ar[r]& \mathcal{L}_{K_{n}(\mathfrak{r})}\ar[r]&
H^{1}_{\mathcal{L}}(K_{n}(\mathfrak{r})_{p},T)\ar[r]&
H^{1}_{\mathcal{L}}(K_{n}(\mathfrak{r})_{p},T)/\mathcal{L}_{K_{n}(\mathfrak{r})}\ar[r]&0}.
$$
 By cohomology we obtain the exact sequence
$$
\xymatrix@=1pc{H^{1}(G_{n,\mathfrak{r}},\mathcal{L}_{K_{n}(\mathfrak{r})})\ar[r]&
H^{1}(G_{n,\mathfrak{r}},H^{1}_{\mathcal{L}}(K_{n}(\mathfrak{r})_{p},T))\ar[r]&
H^{1}(G_{n,\mathfrak{r}},H^{1}_{\mathcal{L}}(K_{n}(\mathfrak{r})_{p},T)/\mathcal{L}_{K_{n}(\mathfrak{r})})\ar[r]&
H^{2}(G_{n,\mathfrak{r}},\mathcal{L}_{K_{n}(\mathfrak{r})})}
$$
Since $\mathcal{L}_{K_{n}(\mathfrak{r})}$ is a free summand of
$\mathbb{V}_{\mathrm{Gal}(\mathcal{K}/K_{n}(\mathfrak{r}))}$ as
$\mathcal{O}[\mathrm{Gal}(K_{n}(\mathfrak{r})/K)]$-module, it
follows that
$$H^{i}(G_{n,\mathfrak{r}},\mathcal{L}_{n}^{\mathfrak{r}})=0\quad \mbox{
for}\;\; i\geq 1. $$
 Hence
$$
H^{1}(G_{n,\mathfrak{r}},H^{1}_{\mathcal{L}}(K_{n}(\mathfrak{r})_{p},T))\cong
H^{1}(G_{n,\mathfrak{r}},H^{1}_{\mathcal{L}}(K_{n}(\mathfrak{r})_{p},T)/\mathcal{L}_{K_{n}(\mathfrak{r})}).
$$
Since
$p^{d}.(H^{1}_{\mathcal{L}}(K_{n}(\mathfrak{r})_{p},T)/\mathcal{L}_{K_{n}(\mathfrak{r})})=0$
(see Lemma \ref{pd s local condition}), we get
$$
p^{d}.H^{1}(G_{n,\mathfrak{r}},H^{1}_{\mathcal{L}}(K_{n}(\mathfrak{r})_{p},T))=0.
\quad \quad \square
$$
Let $c_{n}(\mathfrak{r})\in
\mathbb{Q}_{p}\otimes\bigwedge^{r}H^{1}(K_{n}(\mathfrak{r}),T)$ be
the element defined in $(\ref{definition of Engr})$ and let $
\mathrm{loc}_{p}$ denote the localization map into the semi-local
cohomology at $p$
$$\xymatrix@=2pc{ \mathrm{loc}_{p}:\;
\mathbb{Q}_{p}\otimes
H^{1}(K_{n}(\mathfrak{r}),T)\ar[r]&\mathbb{Q}_{p}\otimes
H^{1}(K_{n}(\mathfrak{r})_{p},T)}.
$$
Since
$\mathbb{Q}_{p}\otimes_{\mathbb{Z}_{p}}\mathbb{V}_{\mathcal{G}_{K_{n}(\mathfrak{r})}}
\cong \mathbb{Q}_{p}\otimes_{\mathbb{Z}_{p}}
H^{1}(K_{n}(\mathfrak{r})_{p},T)$, it follows that
$$
\mathrm{loc}^{(r)}_{p}(c_{n}(\mathfrak{r}))\in
\mathbb{Q}_{p}\otimes_{\mathbb{Z}_{p}}\bigwedge^{r}\mathbb{V}_{\mathcal{G}_{K_{n}(\mathfrak{r})}}.
$$
The defining (integrality) property of the Rubin-Stark elements
shows that, for any
$$
\psi=\psi_{1}\wedge\cdots\wedge \psi_{r}\in
\bigwedge^{r}\mathrm{Hom}(\mathbb{V}_{\mathcal{G}_{K_{n}(\mathfrak{r})}},\mathcal{O}[\mathrm{Gal}(K_{n}(\mathfrak{r})/K)])
$$
we have
$$
\psi(\mathrm{loc}^{(r)}_{p}(c_{n}(\mathfrak{r}))\in
\mathcal{O}[\mathrm{Gal}(K_{n}(\mathfrak{r})/K)].
$$
Hence, by Example $1$ following Proposition $1.2$ in \cite{Rubin96},
we get
$$\mathrm{loc}^{(r)}_{p}(c_{n}(\mathfrak{r}))\in
\bigwedge^{r}\mathbb{V}_{\mathcal{G}_{K_{n}(\mathfrak{r})}}.
$$
Let $F$ be a finite extension of $K$ in $\mathcal{K}$, the map
$$
\xymatrix@=2pc{\mathrm{loc}_{p}:  H^{1}(F,T)\ar[r]& H^{1}(F_{p},T)}.
$$
induces a map
$$
\xymatrix@=1pc{\displaystyle{\varprojlim_{K\subset_{f}F\subset
\mathcal{K}}}\bigwedge^{r-1}\mathrm{Hom}_{\mathcal{O}[\Delta_{F}]}(H^{1}(F_{p},T),
\mathcal{O}[\Delta_{F}])\ar[r]&
\displaystyle{\varprojlim_{K\subset_{f}F\subset
\mathcal{K}}}\bigwedge^{r-1}\mathrm{Hom}_{\mathcal{O}[\Delta_{F}]}(H^{1}(F,T),
\mathcal{O}[\Delta_{F}])}
$$
 The image of
$\Psi\in\displaystyle{\varprojlim_{K\subset_{f}F\subset
\mathcal{K}}}\bigwedge^{r-1}\mathrm{Hom}_{\mathcal{O}[\Delta_{F}]}(H^{1}(F_{p},T),
\mathcal{O}[\Delta_{F}])$ in
$\displaystyle{\varprojlim_{K\subset_{f}F\subset
\mathcal{K}}}\bigwedge^{r-1}\mathrm{Hom}_{\mathcal{O}[\Delta_{F}]}(H^{1}(F,T),
\mathcal{O}[\Delta_{F}])$ will still be denoted by $\Psi$. Let
$\Psi^{\prime}$ be the element constructed in
 Proposition \ref{pro element psi0} and let $\{\varepsilon_{n,\Psi^{\prime}}(\mathfrak{r})^{\chi}\}_{n,\mathfrak{r}}$
 be the Euler system for $T$ associated to $\Psi^{\prime}$  (see Proposition \ref{Euler
 system}). Using Proposition \ref{pro element psi0} and the fact that  $\mathrm{loc}^{(r)}_{p}(c_{n}(\mathfrak{r}))\in
\bigwedge^{r}\mathbb{V}_{\mathcal{G}_{K_{n}(\mathfrak{r})}}$, we see
that
\begin{equation}\label{locp of Euler system and local condition}
\mathrm{loc}_{p}(\varepsilon_{n,\Psi^{\prime}}(\mathfrak{r})^{\chi})\in
H^{1}_{\mathcal{L}}(K_{n}(\mathfrak{r})_{p},T).
\end{equation}
Let $\mathcal{G}$ and $\mathcal{F}$ be Selmer structures on $T$.
Following \cite[\S 2.1]{MR04}, we say that $\mathcal{G}\leq
\mathcal{F}$ if
$$
H^{1}_{\mathcal{G}}(K_{v},T)\subset
H^{1}_{\mathcal{F}}(K_{v},T)\quad\mbox{for all prime $v$.}
$$
 If $\mathcal{G}\leq \mathcal{F}$ we
have an exact sequence \cite[Theorem 2.3.4]{MR04}
$$
\xymatrix@=1.5pc{H^{1}_{\mathcal{G}}(K,T)\ar@{^{(}->}[r]&
H^{1}_{\mathcal{F}}(K,T)\ar[r]&
\bigoplus_{v}H^{1}_{\mathcal{F}}(K_{v},T)/H^{1}_{\mathcal{G}}(K_{v},T)\ar[r]&
H^{1}_{\mathcal{G}^{\ast}}(K,T^{\ast})^{\vee}\ar@{->>}[r]&H^{1}_{\mathcal{F}^{\ast}}
(K,T^{\ast})^{\vee}.}
$$\vskip 6pt
The following lemma is crucial for our purpose
\begin{lem}\label{lemma kolyvagin class}
 Let
 $\kappa_{[K_{n},\mathfrak{r},M]}$ denote the Kolyvagin's derivative
 class, associated to the Euler system $\mathbf{c}=\{p^{d}.\varepsilon_{n,\Psi^{\prime}}(\mathfrak{r})^{\chi}\}_{n,\mathfrak{r}}$, constructed in
 \cite[Chap \textrm{IV},\S 4]{Rubin00}. Then
 $$
\kappa_{[K_{n},\mathfrak{r},M]}\in
H^{1}_{\mathcal{L}^{\mathfrak{r}}}(K_{n},W_{M})
$$
\end{lem}
\noindent \textbf{Proof.} Since $\mathcal{L}^{\mathfrak{r}}\leq
\mathcal{F}_{can}^{\mathfrak{r}}$, we have an exact sequence
$$
\xymatrix@=2pc{
H^{1}_{\mathcal{L}^{\mathfrak{r}}}(K_{n},W_{M})\ar@{^{(}->}[r]&
H^{1}_{\mathcal{F}_{can}^{\mathfrak{r}}}(K_{n},W_{M})\ar[r]&
H^{1}(K_{n,p},W_{M})/H^{1}_{\mathcal{L}}(K_{n,p},W_{M})}.
$$
Theorem $4.5.1$ of \cite{Rubin00} shows that
$\kappa_{[K_{n},\mathfrak{r},M]}\in
H^{1}_{\mathcal{F}_{can}^{\mathfrak{r}}}(K_{n},W_{M})$. Then it
suffices to prove that
$$
\mathrm{loc}_{p}(\kappa_{[K_{n},\mathfrak{r},M]})\in
H^{1}_{\mathcal{L}}(K_{n,p},W_{M})
$$
 where $\mathrm{loc}_{p}$ is the localization map into
the semi-local cohomology at $p$. Let $D_{\mathfrak{r}}$ denote the
derivative operator, defined as in \cite[Definition
\textrm{IV}.4.1]{Rubin00}. Since
 $$
 \xymatrix@=2pc{\mathrm{loc}_{p}:
 H^{1}(K_{n}(\mathfrak{r}),T)\ar[r]&
 H^{1}(K_{n}(\mathfrak{r})_{p},T)}
 $$
is Galois equivariant,
$$
\mathrm{loc}_{p}(D_{\mathfrak{r}}p^{d}.\varepsilon_{n,\Psi}(\mathfrak{r})^{\chi})=
D_{\mathfrak{r}}\mathrm{loc}_{p}(p^{d}.\varepsilon_{n,\Psi}(\mathfrak{r})^{\chi}).
$$
Furthermore, using $(\ref{locp of Euler system and local
condition})$, we get
$$\mathrm{loc}_{p}(p^{d}.\varepsilon_{n,\Psi^{\prime}}(\mathfrak{r})^{\chi})\in
H^{1}_{\mathcal{L}}(K_{n}(\mathfrak{r})_{p},T).$$ On the other hand,
by \cite[Lemma 4.4.2]{Rubin00},
$D_{\mathfrak{r}}p^{d}.\varepsilon_{n,\Psi^{\prime}}(\mathfrak{r})^{\chi}\mod
M$ is fixed by $\mathrm{Gal}(K_{n}(\mathfrak{r})/K_{n})$, then
$$
\mathrm{loc}_{p}(D_{\mathfrak{r}}p^{d}.\varepsilon_{n,\Psi^{\prime}}(\mathfrak{r})^{\chi})\mod
M\in (H^{1}_{\mathcal{L}}(K_{n}(\mathfrak{r})_{p},T)/
MH^{1}_{\mathcal{L}}(K_{n}(\mathfrak{r})_{p},T))^{\mathrm{Gal}(K_{n}(\mathfrak{r})/K_{n})}.
$$
By \cite[Lemma 4.4.13]{Rubin00} and Proposition \ref{Pro Cokernel of
WM and power of p} we have
$$
\mathrm{loc}_{p}(\kappa_{[K_{n},\mathfrak{r},M]})\in
H^{1}_{\mathcal{L}}(K_{n,p},W_{M}).
$$
This finishes the proof of the lemma. \hfill $\square$
\section{Proof of Theorem \ref{AMO2}}
Recall that we proved in Proposition \ref{pro element psi0} the
existence of an element
\begin{equation*}
\Psi^{\prime}=\{\psi_{F}^{\prime}\}_{F}\in
\displaystyle{\varprojlim_{K\subset_{f}F\subset
\mathcal{K}}}\bigwedge^{r-1}\mathrm{Hom}_{\mathcal{O}[\Delta_{F}]}(H^{1}(F_{p},T),\mathcal{O}[\Delta_{F}])
\end{equation*}
such that $\Psi^{\prime}_{F}(\bigwedge^{r}H^{1}(F_{p},T))\subset
H^{1}_{\mathcal{L}}(F_{p},T)$. Let
$\mathbf{c}_{K,\infty}=\{p^{d}.\varepsilon_{n,\Psi^{\prime}}^{\chi}\}_{n}
\in H^{1}_{Iw}(K,T)$ denote the element corresponding to the Euler
system
$\{p^{d}.\varepsilon_{n,\Psi^{\prime}}(\mathfrak{r})^{\chi}\}_{n,\mathfrak{r}}$
 in $H^{1}_{Iw}(K,T):=\displaystyle{\varprojlim_{n}}H^{1}(K_{n},T)$.\\
Remark that under the Leopoldt conjecture for $L$,  the localization
map $$\xymatrix@=2pc{ \mathrm{loc}_{p} : H^{1}(K,T)\ar[r]&
H^{1}(K_{p},T)}$$ is injective. Then by Remark \ref{remark indice
fini} and Proposition \ref{Proposition of locp}, we can find an
element
$\Psi^{\prime}$ such that $\mathbf{c}_{K,\infty}\neq 0$.\\
 One of the keys of the proof of the
main theorem of this article is the following result
\begin{theo}\label{Theo modified condition and char}
Suppose
 the hypotheses $(\mathcal{H}_{0})$ and $(\mathcal{H}_{3})$ hold.
 Then
$$
\mathrm{char}(H^{1}_{\mathcal{L}^{\ast}}(K_{\infty},T^{\ast})^{\vee})\quad\mbox{divides}\quad
\mathrm{char}(H^{1}_{\mathcal{L}}(K_{\infty},T)/\Lambda.\mathbf{c}_{K,\infty})
$$
\end{theo}
\noindent \textbf{Proof.}
 Remark that for any place $v\nmid p$
$$
H^{1}_{\mathcal{L}}(F_{v},T)=H^{1}_{\mathcal{F}_{can}}(F_{v},T).
$$
  Then the proof of this theorem is
similar to the proof of \cite[Theorem 2.3.3]{Rubin00} if we replace
$$
S_{\Sigma_{p}}(F,W^{\ast}_{M})\quad\mbox{by}\quad
H^{1}_{\mathcal{L}}(T,T^{\ast}[M])
$$
and
$$
S^{\Sigma_{p\mathfrak{r}}}(F,W_{M})\quad \mbox{by}\quad
H^{1}_{\mathcal{L}^{r}}(F,W_{M}).
$$
We only need to justify the  following facts:
\begin{enumerate}[label=(\roman*)]
    \item $
\kappa_{[K_{n},\mathfrak{r},M]}\in
H^{1}_{\mathcal{L}^{\mathfrak{r}}}(K_{n},W_{M}).$
    \item $(H^{1}_{\mathcal{L}^{\ast}}(K_{\infty},T^{\ast})^{\vee})_{\Gamma_{n}}\quad
\mbox{and}\quad
\Lambda_{\Gamma_{n}}/\mathrm{char}(H^{1}_{\mathcal{L}^{\ast}}(K_{\infty},T^{\ast})^{\vee})\Lambda_{\Gamma_{n}}
\quad\mbox{are finite}$.
\end{enumerate}
The assertion $(i)$ is Lemma \ref{lemma kolyvagin class}. For the
assertion $(ii)$, on the one hand, we have a surjective map
$$
\xymatrix@=2pc{
H^{1}_{\mathcal{F}_{str}^{\ast}}(K_{\infty},T^{\ast})^{\vee}\ar@{->>}[r]&
H^{1}_{\mathcal{L}^{\ast}}(K_{\infty},T^{\ast})^{\vee}}.
$$
On the other hand, the kernel of the restriction
$$
\xymatrix@=2pc{
H^{1}_{\mathcal{F}_{str}^{\ast}}(K_{\infty},T^{\ast})\ar[r]&
H^{1}_{\mathcal{F}_{str}^{\ast}}(L_{\infty},T^{\ast})}
$$
is finite. Then  it suffices to prove that
$H^{1}_{\mathcal{F}_{str}^{\ast}}(L_{\infty},T^{\ast})^{\Gamma_{n}}$
is finite. Let $M_{\infty}$ be the maximal abelian $p$-extension of
$L_{\infty}$ which is unramified outside the primes above $p$ and
let $\mathfrak{X}_{\infty}=\mathrm{Gal}(M_{\infty}/L_{\infty})$. We
may identify
$H^{1}_{\mathcal{F}_{str}^{\ast}}(L_{\infty},\mathbb{Q}_{p}/\mathbb{Z}_{p})^{\vee}$
with  $\mathfrak{X}_{\infty}$ c.f.\, \cite[\S I.6.3]{Rubin00}.\\
 The Leopoldt conjecture
for $L_{n}$ shows that $\mathfrak{X}_{\infty}^{\Gamma_{n}}=0 $
e.g.\, \cite[Proposition 11.3.3]{NSW91}. Then
$$
H^{1}_{\mathcal{F}_{str}^{\ast}}(L_{\infty},T^{\ast})^{\Gamma_{n}}\quad\mbox{is
finite}.\quad\quad\quad \square
$$
\begin{pro}\label{pro Euler system and strict condition} Suppose
that the hypotheses $(\mathcal{H}_{0})$ and $(\mathcal{H}_{3})$
hold.
 Then

$$
\mathrm{char}(H^{1}_{\mathcal{F}_{str}^{\ast}}(K_{\infty},T^{\ast})^{\vee})\quad\mbox{divides}\quad
\mathrm{char}(H^{1}_{Iw,\mathcal{L}}(K_{p},T)/\mathrm{loc}_{p}(\mathbf{c}_{K,\infty}))
$$
\end{pro}
\noindent \textbf{Proof.} Since $\mathcal{F}_{str}\leq \mathcal{L}$,
we have an exact sequence
$$
\xymatrix@=2pc{
H^{1}_{\mathcal{F}_{str}}(K_{\infty},T)\ar@{^{(}->}[r]&
H^{1}_{\mathcal{L}}(K_{\infty},T)\ar[r]^-{\mathrm{loc}_{p}}&
H^{1}_{Iw,\mathcal{L}}(K_{p},T)\ar[r]&
H^{1}_{\mathcal{F}_{str}^{\ast}}(K_{\infty},T^{\ast})^{\vee}\ar@{->>}[r]&
H^{1}_{\mathcal{L}^{\ast}}(K_{\infty},T^{\ast})^{\vee}}.
$$
 Lemma $\ref{Leopoldt}$ shows that
 $H^{1}_{\mathcal{F}_{str}}(K_{\infty},T)=0$. Then we have an exact sequence
$$
\xymatrix@=1.5pc{
0\ar[r]&H^{1}_{\mathcal{L}}(K_{\infty},T)/\Lambda.\mathbf{c}_{K,\infty}\ar[r]
&H^{1}_{Iw,\mathcal{L}}(K_{p},T)/\mathrm{loc}_{p}(\mathbf{c}_{K,\infty})\ar[r]&
H^{1}_{\mathcal{F}_{str}^{\ast}}(K_{\infty},T^{\ast})^{\vee}\ar@{->>}[r]&
H^{1}_{\mathcal{L}^{\ast}}(K_{\infty},T^{\ast})^{\vee}}.
$$
Theorem $\ref{Theo modified condition and char}$ permits to
conclude.\hfill $\square$\vskip 6pt
Since $K_{\infty}$ is the cyclotomic $\mathbb{Z}_{p}$-extension of
$K$, the  finite primes of $K$ do not split completely in
$K_{\infty}/K$. Therefore, taking the inverse limit in Lemma \ref{pd
s local condition} we deduce that  the $\Lambda$-modules
\begin{equation*}
\displaystyle{\varprojlim_{K\subset F\subset
    K_{\infty}}}\mathcal{L}_{F}\quad\mbox{and}\quad H^{1}_{Iw,\mathcal{L}}(K_{p},T):=\displaystyle{\varprojlim_{K\subset F\subset
    K_{\infty}}}H^{1}_{\mathcal{L}}(F_{p},T)
\end{equation*}
 are pseudo-isomorphic.
  The $\Lambda$-modules
\begin{equation*}
H^{1}_{Iw}(K_{p},T)\quad\mbox{and}\quad
\displaystyle{\varprojlim_{K\subset F\subset
K_{\infty}}}\mathbb{V}_{\mathcal{G}_{F}}
\end{equation*} are also pseudo-isomorphic, thanks to Proposition \ref{pro  coinvariant
and cokernel}. Therefore, by Proposition \ref{pro Psi and
coinvariant}, we conclude that the $\Lambda$-modules
\begin{equation}\label{remark pseudo isomorphic}
H^{1}_{Iw,\mathcal{L}}(K_{p},T)\quad \mbox{and}\quad
\bigwedge^{r}H^{1}_{Iw}(K_{p},T)\;\mbox{are pseudo-isomorphic.}
\end{equation}
 Let $\iota$ denote the composite
of the natural maps
$$
\xymatrix@=2pc{\bigwedge^{r}\displaystyle{\varprojlim_{n}}H^{1}(K_{n},T)\ar[r]&
\displaystyle{\varprojlim_{n}}\bigwedge^{r}H^{1}(K_{n},T)\ar[r]&
\displaystyle{\varprojlim_{n}}(\mathbb{Q}_{p}\otimes_{\mathbb{Z}_{p}}\bigwedge^{r}H^{1}(K_{n},T))}
$$
and let $c_{\infty}:=\{c_{n}\}_{n\geq
0}\in\displaystyle{\varprojlim_{n}}(\mathbb{Q}_{p}\otimes_{\mathbb{Z}_{p}}\bigwedge^{r}H^{1}(K_{n},T))$,
where $c_{n}$ is
 defined in Definition $(\ref{definition of Engr})$.
\begin{theo}\label{theorem fcan divides d Euler system}
 Let $\mathbf{c}$ be an element in
$\iota^{-1}(p^{d}.c_{\infty})$. Under the hypotheses
$(\mathcal{H}_{0})$ and $(\mathcal{H}_{3})$,
$$
\mathrm{char}(H^{1}_{\mathcal{F}_{can}^{\ast}}(K_{\infty},T^{\ast})^{\vee})\quad\mbox{divides}\quad
\mathrm{char}((\bigwedge^{r}H^{1}_{\mathcal{F}_{can}}(K_{\infty},T))/\Lambda.\mathbf{c})
$$
\end{theo}
\noindent \textbf{Proof.} Since $\mathcal{F}_{str}\leq
\mathcal{F}_{can}$ and $H^{1}_{\mathcal{F}_{str}}(K_{\infty},T)=0$
(see Proposition \ref{Leopoldt}), we  have an exact sequence
$$
\xymatrix@=2pc{
H^{1}_{\mathcal{F}_{can}}(K_{\infty},T)\ar@{^{(}->}[r]^-{\mathrm{loc}_{p}}&
H^{1}_{Iw}(K_{p},T)\ar[r]&
H^{1}_{\mathcal{F}_{str}^{\ast}}(K_{\infty},T^{\ast})^{\vee}\ar@{->>}[r]&
H^{1}_{\mathcal{F}_{can}^{\ast}}(K_{\infty},T^{\ast})^{\vee}}.
$$
Proposition \ref{Proposition fcan is free of rank r} (resp. Lemma
\ref{Lemma structur of H(Kp)})) shows that the $\Lambda$-module
$H^{1}_{\mathcal{F}_{can}}(K_{\infty},T)$ (resp.
$H^{1}_{Iw}(K_{p},T)$) is free of rank $r$, then the injection
$\xymatrix@=2pc{
H^{1}_{\mathcal{F}_{can}}(K_{\infty},T)\ar@{^{(}->}[r]^-{\mathrm{loc}_{p}}&
H^{1}_{Iw}(K_{p},T)}$ induces an exact sequence:
$$
\xymatrix@=2pc{ 0\ar[r]&
(\bigwedge^{r}H^{1}_{\mathcal{F}_{can}}(K_{\infty},T))/\Lambda.\mathbf{c}\ar[r]&
(\bigwedge^{r}H^{1}_{Iw}(K_{p},T))/\Lambda.(\mathrm{loc}_{p}^{(r)}(\mathbf{c}))\ar@{->>}[r]&
\mathrm{coker}(\mathrm{loc}^{(r)}_{p})},
$$
where $\mathrm{loc}^{(r)}_{p}$ denotes the map induced on the $r$-th
exterior power. Hence
 $$
 \mathrm{char}((\bigwedge^{r}H^{1}_{Iw}(K_{p},T))/\Lambda.(\mathrm{loc}_{p}^{r}(\mathbf{c})))=
 \mathrm{char}((\bigwedge^{r}H^{1}_{\mathcal{F}_{can}}(K_{\infty},T))/\Lambda.\mathbf{c}).
 \mathrm{char}(\mathrm{coker}(\mathrm{loc}^{(r)}_{p})).
 $$
On the other hand, using $(\ref{remark pseudo isomorphic})$, we see
that the $\Lambda$-modules\\
$\bigwedge^{r}H^{1}_{Iw}(K_{p},T)/\Lambda.\mathrm{loc}_{p}^{(r)}(\mathbf{c}))$
and
$H^{1}_{Iw,\mathcal{L}}(K_{p},T)/\Lambda.\mathrm{loc}_{p}(\mathbf{c}_{K,\infty})$
are pseudo-isomorphic. Since

$$
\mathrm{char}(H^{1}_{\mathcal{F}_{str}^{\ast}}(K_{\infty},T^{\ast})^{\vee})=
\mathrm{char}(H^{1}_{\mathcal{F}_{can}^{\ast}}(K_{\infty},T^{\ast})^{\vee}).
\mathrm{char}(\mathrm{coker}(\mathrm{loc}_{p})),
$$
it follows that
$$
\mathrm{char}(H^{1}_{\mathcal{F}_{str}^{\ast}}(K_{\infty},T^{\ast})^{\vee})\quad\mbox{divides}\quad
\mathrm{char}(H^{1}_{Iw,\mathcal{L}}(K_{p},T)/\Lambda.\mathrm{loc}_{p}(\mathbf{c}_{K,\infty}))
$$
(see Proposition \ref{pro Euler system and strict condition}). Hence
the result follows from the fact that
$$
\mathrm{char}(\mathrm{coker}(\mathrm{loc}^{(r)}_{p}))=\mathrm{char}(\mathrm{coker}(\mathrm{loc}_{p}))
$$
see \cite[page 258]{Bourbaki}. \hfill $\square$\vskip 6pt
\begin{pro}\label{Theorem Assim}
Under the assumption $(\mathcal{H}_{3})$, the cokernel of
$$
\xymatrix@=2pc{
(\widehat{\mathcal{E}}_{\infty}\otimes_{\mathbb{Z}_{p}}\mathcal{O}(\chi^{-1}))_{\Delta}
\ar[r]^-{N_{\Delta}}&(\widehat{\mathcal{E}}_{\infty}\otimes_{\mathbb{Z}_{p}}\mathcal{O}(\chi^{-1}))^{\Delta}
}
$$
is pseudo-null, where $\Delta=\mathrm{Gal}(L_{\infty}/K_{\infty})$.
\end{pro}
\noindent \textbf{Proof.} This is Theorem $5.13$ of
\cite{AMO1}.\hfill $\square$
\begin{coro}\label{coker of power N Delta pseudo nul}
Under the assumption $(\mathcal{H}_{3})$, the cokernel of
$$
\xymatrix@=2pc{
\bigwedge^{r}(\widehat{\mathcal{E}}_{\infty}\otimes_{\mathbb{Z}_{p}}\mathcal{O}(\chi^{-1}))_{\Delta}
\ar[r]^-{N^{(r)}_{\Delta}}&\bigwedge^{r}(\widehat{\mathcal{E}}_{\infty}\otimes_{\mathbb{Z}_{p}}\mathcal{O}(\chi^{-1}))^{\Delta}
}
$$
is pseudo-null.
\end{coro}
\noindent \textbf{Proof.} Let $\mathfrak{p}$ be a prime ideal of
$\Lambda$ of height $\leq 1$. By Proposition \ref{Theorem Assim},
the cokernel of
$$
\xymatrix@=2pc{
(\widehat{\mathcal{E}}_{\infty}\otimes_{\mathbb{Z}_{p}}\mathcal{O}(\chi^{-1}))_{\Delta}
\ar[r]^-{N_{\Delta}}&(\widehat{\mathcal{E}}_{\infty}\otimes_{\mathbb{Z}_{p}}\mathcal{O}(\chi^{-1}))^{\Delta}
}
$$ is pseudo-null, then
$$\mathrm{Im}(N_{\Delta})_{\mathfrak{p}}\cong
((\widehat{\mathcal{E}}_{\infty}\otimes_{\mathbb{Z}_{p}}\mathcal{O}(\chi^{-1}))^{\Delta})_{\mathfrak{p}},
$$
so
$$
(\mathrm{Im}(N^{(r)}_{\Delta}))_{\mathfrak{p}}=\bigwedge^{r}\mathrm{Im}(N_{\Delta})_{\mathfrak{p}}\cong
\bigwedge^{r}((\widehat{\mathcal{E}}_{\infty}\otimes_{\mathbb{Z}_{p}}\mathcal{O}(\chi^{-1}))^{\Delta})_{\mathfrak{p}}.
$$
It follows that the cokernel of
$$
\xymatrix@=2pc{
\bigwedge^{r}(\widehat{\mathcal{E}}_{\infty}\otimes_{\mathbb{Z}_{p}}\mathcal{O}(\chi^{-1}))_{\Delta}
\ar[r]^-{N^{(r)}_{\Delta}}&\bigwedge^{r}(\widehat{\mathcal{E}}_{\infty}\otimes_{\mathbb{Z}_{p}}\mathcal{O}(\chi^{-1}))^{\Delta}
}
$$
is pseudo-null.\hfill $\square$\vskip 6pt
  Recall that $St_{n}$ denotes the
  $\mathbb{Z}[\mathrm{Gal}(L_{n}/K)]$-module generated  by the
  Rubin-Stark elements (see Definition \ref{Definition of Rubin Strak
  module}). Recall also that
$$
c_{n}=\mathrm{cor}^{(r)}_{L_{n+1},K_{n}}(\widetilde{\varepsilon}_{n+1,\chi})
$$
denotes the element defined in $(\ref{definition of Engr})$. Let
$St_{\infty}:=\varprojlim_{n}St_{n}$ and let
$\widetilde{\varepsilon}_{\infty,\chi}:=\{\widetilde{\varepsilon}_{n,\chi}\}_{n\geq
1}$. Since for $n\geq 1$, $
c_{n}=\mathrm{cor}^{(r)}_{L_{n},K_{n}}(\widetilde{\varepsilon}_{n,\chi})
$, it follows that
\begin{eqnarray*}
  \mathrm{res}^{(r)}_{K_{n},L_{n}}(c_{n}) &=& \mathrm{res}^{(r)}_{K_{n},L_{n}}(\mathrm{cor}^{(r)}_{L_{n},K_{n}}(\widetilde{\varepsilon}_{n,\chi})) \\
   &=& |\Delta|^{r-1}N_{\Delta}(\widetilde{\varepsilon}_{n,\chi})
\end{eqnarray*}
 Therefore, using the fact that the restriction map
$$
\xymatrix@=2pc{\mathrm{res}_{K_{n},L_{n}} : H^{1}(K_{n},T)\ar[r]&
H^{1}(L_{n},T)^{\mathrm{Gal}(L_{n}/K_{n})}}
$$ is an isomorphism, see  $(\ref{lemme non semi simple})$. We obtain
$$
|\Delta|^{r-1}N_{\Delta}((\widehat{St_{\infty}})_{\chi})= \Lambda
\mathbf{c},
$$
where $\mathbf{c}$ is the inverse image of
$|\Delta|^{r-1}N_{\Delta}(\widetilde{\varepsilon}_{\infty,\chi})$
under the composite
$$
\xymatrix@=2pc{\bigwedge^{r}\varprojlim_{n}H^{1}(K_{n},T)\ar[r]&
\varprojlim_{n}\bigwedge^{r}H^{1}(K_{n},T)\ar[r]&
\varprojlim_{n}(\mathbb{Q}_{p}\otimes_{\mathbb{Z}_{p}}\bigwedge^{r}H^{1}(K_{n},T))}.
$$
Recall that
$$
H^{1}_{\mathcal{F}_{can}}(K_{\infty},T)\cong
(\widehat{\mathcal{E}^{\prime}_{\infty}}\otimes_{\mathbb{Z}_{p}}\mathcal{O}(\chi^{-1}))^{\Delta}
\quad(\mbox{see}\; (\ref{Fcan and units })).
$$
 \noindent \textbf{Proof Theorem \ref{AMO2}.}
 Consider the commutative exact diagram
$$
\xymatrix@=1.5pc{ &
(\widehat{\mathrm{St}_{\infty}})_{\chi}\ar[r]\ar@{->>}[d]^-{|\Delta|^{r-1}N_{\Delta}}&
\bigwedge^{r}(\widehat{\mathcal{E}_{\infty}})_{\chi}\ar@{->>}[r]\ar[d]^-{N^{(r)}_{\Delta}}&
\big(\bigwedge^{r}\widehat{\mathcal{E}_{\infty}}/\widehat{\mathrm{St}_{\infty}}\big)_{\chi}\ar[d]\\
0\ar[r]&|\Delta|^{r-1}N_{\Delta}(\mathbf{c})\ar[r]&
\bigwedge^{r}(\widehat{\mathcal{E}_{\infty}})^{\chi}\ar@{->>}[r]\ar@{->>}[d]&
\bigwedge^{r}(\widehat{\mathcal{E}_{\infty}})^{\chi}/|\Delta|^{r-1}N_{\Delta}(\mathbf{c})\\
  &  & \mathrm{coker}(N^{(r)}_{\Delta})& }
$$
where
$(\widehat{\mathcal{E}_{\infty}})^{\chi}=(\widehat{\mathcal{E}_{\infty}}\otimes_{\mathbb{Z}_{p}}
\mathcal{O}(\chi^{-1}))^{\Delta}$. Corollary \ref{coker of power N
Delta pseudo nul}
 shows that the $\Lambda$-module $\mathrm{coker}(N^{(r)}_{\Delta})$
 is pseudo-null, so
 $$
 \mathrm{char}(\bigwedge^{r}(\widehat{\mathcal{E}_{\infty}})^{\chi}/
 |\Delta|^{r-1}N_{\Delta}(\mathbf{c}))
\quad\mbox{divides}\quad
\mathrm{char}\bigg(\big(\bigwedge^{r}\widehat{\mathcal{E}_{\infty}}/\widehat{\mathrm{St}_{\infty}}\big)_{\chi}\bigg).
$$
Since $\chi(D_{v}(L/K))\neq 1$ for any $p$-adic prime of $K$, then
$$
(\widehat{\mathcal{E}^{\prime}_{\infty}}\otimes_{\mathbb{Z}_{p}}\mathcal{O}(\chi^{-1}))^{\Delta}
\cong
(\widehat{\mathcal{E}_{\infty}}\otimes_{\mathbb{Z}_{p}}\mathcal{O}(\chi^{-1}))^{\Delta}.
$$
Hence the theorem follows from Theorem \ref{theorem fcan divides d
Euler system} and Lemma \ref{chi quotient of class group and
canonique conditio};
$$
\mathrm{char}((A_{\infty})_{\chi})\quad \mbox{divides}\quad p^{d}.
\mathrm{char}
\bigg(\big(\bigwedge^{r}\widehat{\mathcal{E}_{\infty}}/\widehat{\mathrm{St}_{\infty}}\big)_{\chi}\bigg).
$$
\vskip 10pt
\begin{tabbing}
  \hspace{6cm}\=\hspace{3cm} \= \hspace{5cm} \= \kill
Département de mathématiques\>  \> Laboratoire de mathématique\\
Faculté des sciences de Meknes\> \> 16 Route de Gray\\
B.P. 11201 Zitoune\>  \> 25030 Besançon\\
Meknes, Maroc\>  \> cedex, France\\
 mazigh.younes@gmail.com\>\>youness.mazigh@univ-fcomte.fr
  \end{tabbing}

\begin{thebibliography}{99}
     \setlength{\parskip}{0cm}
  \bibitem[AMO]{AMO1}{\bf Assim,\,J. Mazigh, Y. Oukhaba, H.} {\it  Théorie d'Iwasawa des unités de Stark et groupe de
  classes.} to appear in  International Journal of Number Theory.
  \bibitem[Bo]{Bourbaki} {\bf Bourbaki. N} {\it Algèbre commutative, chapitre $7$, Diviseurs}, Hermann, 1965
   \bibitem[B\"{u} 07]{Kazim107}{\bf B\"{u}y\"{u}kboduk, K.} {\it $\Lambda$-adic Kolyvagin systems}. Int. Math. Res. Not. IMRN 2011, no. 14, 3141-3206
 \bibitem[B\"{u} 08]{Kazim108}{\bf B\"{u}y\"{u}kboduk, K.} {\it Kolyvagin systems of Stark
 units }. J. Reine Angew. Math. 631 (2009), 85-107.
     \bibitem[B\"{u} 09]{Kazim109}{\bf B\"{u}y\"{u}kboduk, K.} {\it Stark units and the main conjectures for totally real fields}.
      Compos. Math. 145 (2009), no. 5, 1163-1195.

 \bibitem[Gr 94]{Greither94}{\bf  Greither, C.} {\it Sur les normes universelles dans
les $\mathbb{Z}_{p}$-extensions.} Journal de Théorie des Nombres de
Bordeaux 6 (1994), 205-220.
 \bibitem[Gr 96]{Greither96}{\bf  Greither, C.} {\it On Chinburg's second conjecture for abelian fields}
 J. Reine Angew. Math. 479 (1996), 1-37.
\bibitem[Mi 86]{Milne}{\bf Milne, J.}{\it Arithmetic Duality theorems.}
Acad. Press, Boston, 1986.
 \bibitem[MR 04]{MR04}{\bf Mazur, B.; Rubin, K.} {\it Kolyvagin systems}. Mem. Amer. Math. Soc., 168(799):\textrm{viii}+96, 2004.
\bibitem[MR 16 ]{MR 16}{\bf Mazur, B.; Rubin, K.} {\it Controlling Selmer groups in the higher core rank
case}. J. Théor. Nombres Bordeaux 28 (2016), no. 1, 145-183.
\bibitem[NSW 91]{NSW91}{\bf Neukirch,\,J.; Schmidt, A.; Wingberg, K.} {\it Cohomology of Number Fields}. Springer (1991).
\bibitem[Ne 06]{Nekovar06}{\bf Nekov\'{a}\~{r},\,J.} {\it Selmer complexes}. Astérisque 310 (2006).
\bibitem[PR 92]{PR92}{\bf B. Perrin-Riou.} {\it Théorie d'Iwasawa et hauteurs p-adiques}. Invent. Math. 109 (1992) 137-185.
\bibitem[Po 02]{Popescu02}{\bf Popescu, C.} {\it Base change for Stark-type conjectures "over
$\mathbb{Z}$"}. J. Reine Angew. Math, 524:85-111, 2002.
\bibitem[Ru 96]{Rubin96}{\bf  Rubin, K.} {\it  A Stark conjecture "over $\mathbb{Z}$" for abelian $L$-functions with multiple
zeros}. Ann. Inst. Fourier (Grenoble), 46(1):33-62, 1996.
 \bibitem[Ru 00]{Rubin00}  {\bf    Rubin, K.} {\it  Euler systems}. Annals of Mathematics Studies, 147.
        Hermann Weyl Lectures. The Institute for Advanced Study. Princeton
        University Press, Princeton,
 \bibitem[Ta 84]{Tate84} {\bf Tate, J.} {\it Les conjectures de Stark sur les fonctions L d'Artin en s=0}.
Birkh\"{a}user Boston Inc, 1984. Lecture notes edited by Dominique
Bernardi and Norbert Schappacher.
\end{thebibliography}
\end{document}